\documentclass[a4paper,11pt,leqno,english]{amsart}
\usepackage[utf8]{inputenc}
\usepackage[T1]{fontenc}
\usepackage{microtype}
\usepackage[french, main=english]{babel}
\usepackage{amssymb,latexsym}
\newcommand{\mmod}[1]{\ (\mathrm{mod}\ #1)}
\usepackage{mathrsfs,tgschola}
\usepackage[normalem]{ulem}
\usepackage{pgfplots}
\usepackage{tkz-fct}
\usepackage{url}
\usepackage{xcolor}
\usepackage{comment}
\definecolor{violet}{rgb}{0.0,0.2,0.7}
\definecolor{rouge2}{rgb}{0.8,0.0,0.2}
\usepackage{tikz}
\usepackage{empheq}
\usepackage{tikz-cd}
\usetikzlibrary{matrix,arrows,decorations.pathmorphing}
\usepackage{hyperref}
\usepackage{mathpazo}
\usepackage{enumerate}
\usepackage[top=1in, bottom=1in, left=1.13in, right=1.13in]{geometry}

\hypersetup{
    bookmarks=true,         
    unicode=false,          
    pdftoolbar=true,        
    pdfmenubar=true,        
    pdffitwindow=false,     
    pdfstartview={FitH},    
    pdftitle={},    
    pdfauthor={},     
    colorlinks=true,       
   linkcolor=rouge2,          
    citecolor=violet,        
    filecolor=black,      
    urlcolor=cyan}           
\usepackage{enumitem}
\usepackage{appendix}
\usepackage{amsmath}

\makeatletter
\newenvironment{smallalign}{%
  \tiny  
  \let\oldmaketag@@@\maketag@@@
  \def\maketag@@@##1{\hbox{\m@th\normalsize\normalfont##1}}
  \start@align\@ne\st@rredfalse\m@ne
}{%
  \endalign
}
\makeatother

\theoremstyle{definition}
\newtheorem{defi}{Definition}[section]

\newtheorem{nota}[defi]{Remark}

\newtheorem{conj}[defi]{Conjecture}

\theoremstyle{plain}
\newtheorem{thm}[defi]{Theorem}
\newtheorem{prp}[defi]{Proposition}
\newtheorem{lem}[defi]{Lemma}
\newtheorem{coro}[defi]{Corollary}
\newtheorem{spthmA}{Theorem}

\newtheorem{spthmB}{Theorem}

\newcommand{\scr}{\mathscr}
\newcommand{\bb}{\mathbb}
\newcommand{\mcal}{\mathcal}
\newcommand{\spec}{\operatorname{Spec}}
\newcommand{\Hom}{\operatorname{Hom}}
\newcommand{\Ext}{\operatorname{Ext}}

\newcommand{\RR}{\mathbf{R}}

\newcommand{\vir}{\mathsf{vir}}
\newcommand{\Hilb}{\mathsf{Hilb}}
\newcommand{\dHilb}{\mathsf{dHilb}}
\newcommand{\MSh}{\mathsf{M}}
\newcommand{\dMSh}{\mathsf{dM}}

\newcommand{\Sym}{\operatorname{Sym}}
\newcommand{\CH}{\operatorname{CH}}
\newcommand{\vdim}{\operatorname{vdim}}
\newcommand{\fix}{\mathsf{fix}}
\newcommand{\mov}{\mathsf{mov}}
\newcommand{\ncHilb}{\mathsf{ncHilb}}
\newcommand{\ncM}{\mathsf{ncM}}

\newcommand{\DTC}{\mathsf{DT}}
\newcommand{\DT}{\mathsf{DT}}
\newcommand{\Tot}{\operatorname{Tot}}

\newcommand{\Res}{\mathsf{Res}}
\newcommand{\End}{\operatorname{End}}
\newcommand{\GL}{\operatorname{GL}}

\newcommand{\tr}{\operatorname{tr}}
\newcommand{\rank}{\operatorname{rank}}
\newcommand{\coker}{\operatorname{coker}}

\title[Donaldson-Thomas invariants of {$[\mathbb C^4/ {\bb Z_{\MakeLowercase{r}}}]$}]{Donaldson-Thomas invariants of $[\mathbb C^4/\bb Z_r]$}

\author{Xiaolong Liu}
\address{Institute of Mathematics, AMSS, Chinese Academy of Sciences, 55 Zhongguancun East Road, Beijing, 100190, China}
\email{liuxiaolong@amss.ac.cn}
\date{\today}

\begin{document}

\begin{abstract}
We compute the zero-dimensional Donaldson-Thomas invariants of the quotient stack 
$[\bb{C}^4/\bb{Z}_r]$, confirming a conjecture of Cao-Kool-Monavari. Our main theorem is 
established through an orbifold analogue of Cao-Zhao-Zhou's degeneration formula combined 
with the zero-dimensional Donaldson-Thomas invariants for $\widetilde{\bb C^2/\bb Z_r}\times\bb{C}^2$ and 
an explicit determination of orientations of Hilbert schemes of points on $[\bb{C}^4/\bb{Z}_r]$.
\end{abstract}

\maketitle

\hypersetup{linkcolor=black}

\section{Introduction}

\subsection{Background}
Donaldson-Thomas theory is a sheaf-theoretic framework in enumerative geometry, initially developed for Calabi-Yau and Fano $3$-folds, which first appeared in Thomas' thesis \cite{Thomas2000}. The theory is based on the two-term perfect obstruction theories discovered by Li-Tian \cite{LT98} and Behrend-Fantechi \cite{BF97}. Later, Maulik-Nekrasov-Okounkov-Pandharipande \cite{MNOP1,MNOP2} extended the definitions of Donaldson-Thomas invariants for Hilbert schemes to any $3$-folds.
However, unlike Gromov-Witten invariants, extending Donaldson-Thomas theory to higher-dimensional algebraic
varieties presents difficulties. Since obstruction theories on moduli spaces of sheaves are no longer two-term,
the method of constructing virtual cycles in \cite{LT98,BF97} does not work for higher-dimensional varieties.

Donaldson-Thomas theory for Calabi-Yau $4$-folds (or $\mathsf{DT}_4$ theory) was first studied by Cao-Leung \cite{CL14}, who constructed $\mathsf{DT}_4$ virtual classes in several special cases. The first complete mathematical formulation of the theory was established by Borisov-Joyce \cite{BJ17}, based on the derived differential geometry and the shifted symplectic geometry developed by Pantev-Toën-Vaquié-Vezzosi \cite{PTVV13}. Subsequently, an algebraic approach was provided by Oh-Thomas \cite{OT23}, with connections between the two frameworks explored further in \cite{OT24}. More recently, Park \cite{park21,park24} developed virtual pullback theories for $(-2)$-shifted symplectic fibrations, making $\mathsf{DT}_4$ virtual classes functorial.

These fundamental developments have influenced various directions in enumerative geometry, including Donaldson-Thomas invariants appearing in \cite{CK18,CK20,CKM22,monavari22,CKM23,CZZ24,Bojko24} and \cite{KR}; connections to Gromov-Witten, Pandharipande-Thomas and Gopakumar-Vafa invariants in \cite{CMT18,CMT21}; surface counting in \cite{BKP22,BKP24}; and gauged linear sigma model (GLSM) in \cite{CZ23,CTZ25,KP25}.

\subsection{Zero-dimensional \texorpdfstring{$\mathsf{DT}_4$}{DT4}-invariants of toric Calabi-Yau \texorpdfstring{$4$}{4}-folds}
Analogous to the zero-dimensional Donaldson-Thomas invariants for $3$-folds studied in \cite{MNOP1,MNOP2,Li06,LP09}, the $\mathsf{DT}_4$-invariants of toric Calabi-Yau $4$-folds\footnote{We call $X$ a \textit{toric Calabi-Yau $4$-fold} if it is a smooth quasi-projective toric $4$-fold satisfying $K_X\cong\scr O_X$, $H^{>0}(X,\scr O_X)=0$ and every cone of its fan is contained in a $4$-dimensional cone.} are conjectured as follows.

\begin{defi}[$\DT_4$-invariants of toric CY $4$-folds]\label{def-DT4-toric-CY4}
Let $X$ be a toric Calabi-Yau $4$-fold with a $\bb T$-equivariant line bundle $\scr L$ on $X$
where $\bb T\subset(\bb C^*)^4$ is the Calabi-Yau torus.
We define the \textit{tautological bundle of $\scr L$} on $\Hilb^n(X)$ as
\[\scr L^{[n]}:=\RR\pi_*(\pi_X^*\scr L\otimes\scr O_{\mcal Z}),\]
where $\mcal Z\subset X\times\Hilb^n(X)$ is the universal closed subscheme and $\pi:X\times\Hilb^n(X)\to\Hilb^n(X)$,
$\pi_X:X\times\Hilb^n(X)\to X$ are projections. Let $\bb C^*_m$ be an $1$-dimensional torus with trivial action on $X$ and let $y=e^m$ be a character of $\bb C_m^*$ with weight $m=c_1(y)$. Fix an orientation, we define the generating function of \textit{zero-dimensional $\DT_4$-invariants with tautological insertions} as 
\begin{align*}
\DTC_{X,\scr L}(m,q)&:=1+\sum_{n\geq1}q^n\int_{[\Hilb^n(X)]^{\vir}_{\bb T}}e_{\bb T\times\bb C^*_m}((\scr L^{[n]})^{\vee}\otimes e^m)\\
    &\in\operatorname{Frac}\left(\frac{\bb C[s_1,s_2,s_3,s_4]}{(s_1+s_2+s_3+s_4)}\right)[m][[q]]\cong\bb C(s_1,s_2,s_3)[m][[q]],
\end{align*}
where $[\Hilb^n(X)]^{\vir}_{\bb T}$ be the $\bb T$-equivariant virtual class and $s_i=c_1(t_i)$ be the weights of the coordinates of $\bb T$.
\end{defi}

\begin{conj}[{\cite[Conj.~1.6]{CK18}}]\label{conj-torus-DT4-general}
For any toric Calabi-Yau $4$-fold $X$, there exists a choice of orientation such that
\begin{equation*}
    \DTC_{X,\scr L}(m,q)=M(-q)^{\int_Xc_1^{\bb T\times\bb C_m^*}(\scr L\otimes y^{-1})c_3^{\bb T}(X)},
\end{equation*}
where $M(q):=\prod_{n\geq1}(1-q^n)^{-n}$ be the MacMahon series and $\int_X$ be the $\bb T$-equivariant pushforward
to a point.
\end{conj}

The case of $X=\bb C^4$ and the CY local curves was proved by Cao-Zhao-Zhou \cite[Thm.~6.8,~6.16]{CZZ24}.
Note that upon ignoring the orientation issue, this conjecture can be deduced by the vertex calculation 
in \cite{CK18,NP2019} and the result for $\bb C^4$. In this paper, we prove a compatibility of the 
orientations (Proposition \ref{thm-sign-rule-toric}) and deduce the following result.

\begin{spthmA}[Theorem \ref{coro-DT-specil-toric}, Corollary \ref{coro-DT-ArC2}]\label{thm-coro-DT-special-toric}
Let $X=Y\times\bb C$ for some
toric Calabi-Yau $3$-fold $Y$, then there exists a choice of orientation such that Conjecture \ref{conj-torus-DT4-general} holds.

In particular, let $X=\widetilde{\bb C^2/\bb Z_r}\times\bb C^2$
where $\widetilde{\bb C^2/\bb Z_r}$ be the minimal crepant resolution of $A_{r-1}$-singularity 
$\bb C^2/_{\text{aff}}\bb Z_r$. Then there exists a choice of orientation such that
\begin{equation*}
    \DTC_{\widetilde{\bb C^2/\bb Z_r}\times\bb C^2,\scr O}(m,q)=M(-q)^{-\frac{m}{s_4}
    \left(\frac{r(s_1+s_2)}{s_3}+\frac{(s_1+s_2)(s_1+s_2+s_3)}{rs_1s_2}\right)}.
\end{equation*}
\end{spthmA}

\subsection{\texorpdfstring{$\mathsf{DT}_4$}{DT4}-invariants of \texorpdfstring{$[\mathbb C^4/\bb Z_r]$}{C4Zr} and the crepant resolution conjecture}
In \cite{CKM23}, the authors investigated the Donaldson-Thomas theory for toric Calabi-Yau $4$-orbifolds, with $[\bb C^4/G]$ serving as the local model for some finite subgroup $G<\text{SU}(4)$.
Here we consider $G=\bb Z_r<\text{SU}(2)<\text{SU}(4)$ with the following conjecture due to Cao-Kool-Monavari.

\begin{conj}[{\cite[Conj.~5.13,~6.6]{CKM23}}]\label{conj-DT-C4Zr}
There exists a choice of orientation such that
\begin{equation*}
\begin{aligned}
    \DTC_{[\bb C^4/\bb Z_r],\scr O}(m,q_0,...,q_{r-1})
    =&M(1,-Q)^{-\frac{m}{s_4}\left(\frac{r(s_1+s_2)}{s_3}+\frac{(s_1+s_2)(s_1+s_2+s_3)}{rs_1s_2}\right)}\\
    &\cdot\prod_{0<i\leq j<r}\widetilde{M}(q_{[i,j]},- Q)^{-\frac{m(s_1+s_2)}{s_3s_4}},
\end{aligned}
\end{equation*}
where $\DTC_{[\bb C^4/\bb Z_r],\scr O}$ is defined similar as the toric case (see Definition \ref{def-DT4-invariants-C4Zr}), $M(a,b):=\prod_{n=1}^{\infty}(1-ab^n)^{-n}$ is the refined MacMahon series
and $\widetilde{M}(a,b):=M(a,b)\cdot M(a^{-1},b)$.
Here we denote $q_{[i,j]}:=q_iq_{i+1}\cdots q_j$ for $i\leq j$ and $Q=q_{[0,r-1]}$.
\end{conj}
Our main result in this paper is to confirm this conjecture.
\begin{spthmB}[Theorem \ref{thm-main-DT4-C4Zr}]\label{thm-intro-DT4-C4Zr}
Conjecture \ref{conj-DT-C4Zr} holds.
\end{spthmB}

We first compute the sign rules for two canonical orientations in $\S$\ref{sect on ori of C4/Zr} (see Theorem \ref{thm-ori-A-sign-rule} and Corollary \ref{coro-sign-rule-oB}), 
following methods analogous to those of Kool-Rennemo \cite{KR} together with a comparison arguments.  

Next, we construct orbifold analogues of expanded and rubber geometries extending the techniques 
of \cite{Zhou18-2,Zhou18,CZZ24} to carry out several pole analyses in 
\S\ref{sect on DT of C4/Zr}. 
In \S\ref{subsection-toric-cpts}, we compare the absolute $\mathsf{DT}_4$-invariants
of canonical-bundle total spaces over two proper toric surface compactifications.
Together with the preceding relative and rubber pole estimates, this gives further
pole analyses for our invariants.
Finally, our Theorem \ref{thm-intro-DT4-C4Zr} follows from combining these results with 
Theorem \ref{thm-coro-DT-special-toric}; see \S\ref{sec-DT4-C4Zr-final} for details.

In \cite[Conj.~5.16,~6.7]{CKM23}, Cao-Kool-Monavari further conjectured a precise correspondence 
between the $\mathsf{DT}_4$-invariants of the following two resolutions of the affine GIT 
quotient $\bb C^4/_{\text{aff}}\bb Z_r=\spec\bb C[x_1,...,x_4]^{\bb Z_r}$:
\[\begin{tikzcd}
	{[\bb C^4/\bb Z_r]} & {\bb C^4/_{\text{aff}}\bb Z_r} & {\widetilde{\bb C^2/\bb Z_r}\times\bb C^2.}
	\arrow[from=1-1, to=1-2]
	\arrow[from=1-3, to=1-2]
\end{tikzcd}\]
The discrepancy between these invariants is precisely captured by the Pandharipande-Thomas invariants $\mathsf{PT}_{\widetilde{\bb C^2/\bb Z_r}\times\bb C^2,\scr O}$ of $\widetilde{\bb C^2/\bb Z_r}\times\bb C^2$; see \cite[\S 6.1]{CKM23} for details.
This is referred to as the \textit{Donaldson-Thomas crepant resolution conjecture}.

Since the invariants $\DTC_{\widetilde{\bb C^2/\bb Z_r}\times\bb C^2,\scr O}$ and $\DTC_{[\bb C^4/\bb Z_r],\scr O}$ have been computed in 
Theorem \ref{thm-coro-DT-special-toric} and Theorem \ref{thm-intro-DT4-C4Zr}, respectively, the solution to the Donaldson-Thomas crepant 
resolution conjecture reduces to the computation of 
$\mathsf{PT}_{\widetilde{\bb C^2/\bb Z_r}\times\bb C^2,\scr O}$ which we hope to address in a future work.

\subsection{Related works and some future directions}
Kool and Rennemo established the K-theoretic $\mathsf{DT}_4$-invariants for $\mathsf{Quot}^n_r(\bb C^4)$ in their work \cite{KR}, resolving a conjecture posed by Nekrasov \cite{Nekrasov2020} through the application of factorizable sequences of sheaves developed by Okounkov \cite{Okounkov2017}. 
Hence, it is worth considering using a similar approach to compute the K-theoretic 
$\mathsf{DT}_4$-invariants of $[\bb C^4/\bb Z_r]$ conjectured in \cite[Conj.~5.13]{CKM23}; 
see also \cite{Thi24} for the case of $[\bb C^3/\bb Z_r]$. See \cite{AKL26} for related works.

Moreover, the physical paper \cite{ST23} conjectured the higher-rank version, that is, $\mathsf{DT}_4$-invariants of $\mathsf{Quot}([\bb C^4/\bb Z_r])$, which will 
follow from our result and the similar method to \cite[~\S 6]{FMR21}; see also the arguments 
in \cite[~\S 5.3]{KR}.

We observe that similar conjectural expressions exist for $[\bb C^4/(\bb Z_2\times\bb Z_2)]$ in \cite{CKM23}, analogous to Conjecture \ref{conj-DT-C4Zr}. Note that the zero-dimensional $\mathsf{DT}_4$-invariant of the crepant resolution of $\bb C^4/_{\text{aff}}(\bb Z_2\times\bb Z_2)$ could be 
deduced by Theorem \ref{thm-coro-DT-special-toric} as a special case. Therefore, it is worthwhile 
to compute the $\mathsf{DT}_4$-invariants of $[\bb C^4/(\bb Z_2\times\bb Z_2)]$.

Finally, Schmiermann also studied the sign rules of the Hilbert scheme of points on $[\bb C^4/G]$ for the finite subgroup $G\subset\bb T$ in his thesis \cite{Sch25}.

\subsection*{Conventions and Notations}
\begin{enumerate}
    \item All schemes and stacks (classical or derived) are defined over $\bb C$. Chow groups $\CH_*(-)$ and $K_0$-groups $K_0(-)$ are taken with $\bb Q$-coefficients.
    \item The \textit{standard action} of $(\bb C^*)^4$ on $\bb C^4=\spec\bb C[x_1,x_2,x_3,x_4]$ is given by $t\cdot x_i=t_ix_i$ where $t=(t_1,t_2,t_3,t_4)$. At the level of closed points, the action is given by
\[(t_1,t_2,t_3,t_4)\cdot(a,b,c,d)=(t_1^{-1}a,t_2^{-1}b,t_3^{-1}c,t_4^{-1}d).\]
Let $s_i=c_1(t_i)$ be the corrsponding weights.
\item For a space $\mcal{M}$ with a $G$-action, $\mcal{M}^G$ denotes its fixed locus. For a $G$-module $M$, $M^{G\text{-}\fix}$ represents the fixed part. When $G =\bb{T}$ is a torus, $N^{\mathbb{T}\text{-}\mov}$ denotes the moving part of a $\mathbb{T}$-module $N$.
\item For a morphism of stacks (classical or derived) $f:\mcal X\to\mcal Y$, let $\bb L_{f}$ be its full cotangent complex and $\mathbf L_{f}:=\tau^{\geq-1}\bb L_{f}$ be the truncated cotangent complex. Let $\bb T_f:=\bb L_f^{\vee}$ be the tangent complex for 
homotopical smooth morphism $f$.
\end{enumerate}

\subsection*{Acknowledgments}
I am deeply grateful to Yalong Cao for introducing me to this subject and for his patient guidance and many helpful discussions throughout this work. The author also thanks Martijn Kool, Sergej Monavari, and Zijun Zhou for their valuable discussions and comments on the draft of this work.

\section{Recollection of square-root virtual pullbacks}\label{sect on vir pullback}
In this section, we review the foundations of square-root virtual pullback through the 
symmetric obstruction theory developed in \cite{park21}
and its derived enhancement \cite{park24}, which will be used frequently later.

\begin{defi}[{\cite[Prop.~1.7,~\S A.2]{park21}}]\label{defi on sym cx}
A \textit{symmetric complex} $(\bb E,\theta,o)$ on an Artin stack
$\mcal X$ is a perfect complex $\bb E$ of amplitude $[-2,0]$ consists of the following data:
\begin{enumerate}
    \item A non-degenerate symmetric form $\theta$ on $\bb E$, i.e. a morphism
\[\theta : \scr O_{\mcal X} \to (\bb E \otimes\bb E)[-2]\]
in $\mathsf D_{\text{coh}}(\mcal X)$, invariant under the transposition
such that the induced morphism
$\iota_{\theta} : \bb E^\vee\to\bb E[-2]$ is an isomorphism.
    \item An \textit{orientation} $o$ of $\bb E$, i.e., an isomorphism $o:\scr O_\mcal X \xrightarrow{\cong} \det(\bb E)$ such that $\det(\iota_\theta)=o\circ o^\vee$.
\end{enumerate}
\end{defi}

\begin{defi}(\cite[Def.~1.9,~\S A.2]{park21})\label{def on sym ob}
A \textit{symmetric obstruction theory} for a Deligne-Mumford morphism $f : \mcal X\to \mcal Y$ of Artin stacks is a morphism $\phi : \bb E \to \mathbf L_f$ in
$\mathsf D^{\leq0}_{\text{coh}}(\mcal X)$ such that
\begin{enumerate}
    \item $\bb E$ is a symmetric complex defined as in Definition \ref{defi on sym cx}.
    \item $\phi$ is an obstruction theory in the sense of Behrend-Fantechi \cite{BF97}, that is, $h^0(\phi)$ is an isomorphism and $h^{-1}(\phi)$ is surjective.
\end{enumerate}
\end{defi}

For any symmetric complex $\bb E$, there is a quadratic function 
\begin{equation*}
\mathfrak q_\bb E:\mathfrak C_{\bb E}\to \bb A^1_{\mcal X}, 
\end{equation*}
constructed in \cite[Prop.~1.7, \S A.2]{park21},
where $\mathfrak C_{\bb E}:=h^1/h^0((\tau^{\geq-1}\bb E)^{\vee}_{\text{fppf}})$ is the virtual normal cone of $\bb E$, with some natural conditions. 
When $\bb E=E[1]$ for an $\text{SO}(*,\bb C)$-bundle $E$, the function $\mathfrak q_{\bb E}$ is given by the quadratic form on $E$. 

\begin{defi}\label{def on iso sym ob}
A symmetric obstruction theory $\phi : \bb E \to \mathbf L_f$  is called {\it isotropic} if the intrinsic normal
cone $\mathfrak C_f$ is isotropic in the virtual normal cone $\mathfrak C_{\bb E}$, i.e., the restriction
\[\mathfrak q_{\bb E}|_{\mathfrak C_f}
:\mathfrak C_f\hookrightarrow \mathfrak C_{\bb E}\to \bb A^1_{\mcal X}\]
vanishes, where the closed
embedding of the intrinsic normal cone $\mathfrak C_f\hookrightarrow \mathfrak C_{\bb E}$ induced by the obstruction theory.
\end{defi}

\begin{defi}\label{def of virtual pullbacks}
Let $f:\mcal X\to\mcal Y$ be a Deligne-Mumford morphism between Artin stacks with an isotropic symmetric obstruction theory $\phi:\bb E\to\mathbf L_f$.
Let $\mcal X$ be a quotient of a separated Deligne-Mumford stack by an algebraic group.
Let $\mathfrak Q(\bb E)$ be the zero locus of the quadratic function $\mathfrak q_{\bb E}:\mathfrak C_\bb E\to \bb A^1_\mcal X$,
there is a \textit{square-root Gysin pullback} 
$$\sqrt{0^!_{\mathfrak Q(\bb E)}}:\CH_*(\mathfrak Q(\bb E)) \to \CH_*(\mcal X), $$ defined in \cite[Def.~A.2]{park21}.
Since $\phi$ induces a closed embedding $a:\mathfrak C_f\hookrightarrow\mathfrak Q(\bb E)$, the
\textit{square-root virtual pullback} can be defined as the composition 
\begin{equation*}
    \sqrt{f^!}:\CH_*(\mcal Y)\xrightarrow{\text{sp}_f}\CH_*(\mathfrak C_f )\xrightarrow{a_*}\CH_{*} (\mathfrak Q(\bb E))\xrightarrow{\sqrt{0^!_{\mathfrak Q(\bb E)}}}\CH_{*}(\mcal X),
\end{equation*}
where $\text{sp}_f:\CH_*(\mcal Y)\to\CH_*(\mathfrak C_f)$
is the specialization map (\cite[Const.~3.6]{Man12}).

When $\mcal X$ is a separated Deligne-Mumford stack over $\mathbb C$ with the structure morphism $a_{\mcal X}:\mcal X\to\spec\bb C$
and an isotropic symmetric obstruction theory $\phi:\bb E\to\mathbf L_{\mcal X}$, then we define 
the \textit{virtual class} of $X$ as \[[\mcal X]^{\vir}:=\sqrt{a_{\mcal X}^!}[\spec\bb C]\in\CH_{\vdim(\mcal X)}(\mcal X).\]
\end{defi}

The map $\sqrt{f^!}$ commutes with \textit{projective pushforwards, smooth pullbacks, and 
Gysin pullbacks for regular immersions} and can be used to define the square root virtual 
pullback for any base change of $f$. 
Moreover, it is \textit{functorial} with respect to morphisms compatible with symmetric 
obstruction theories; see \cite[\S 2, Thm.~A.4]{park21} for details.

This theory serves as the classical shadow of the virtual Lagrangian pullbacks of locked 
$(-2)$-shifted symplectic fibrations developed in \cite[Thm.~5.2.2]{park24}. Indeed, let $g:\mathsf{M}\to\mathsf{B}$ be an oriented $\mathsf w$-locked $(-2)$-shifted symplectic 
fibration of derived stacks with similar assumptions as before and underlying function $\mathsf w:\mathsf{B}\to\bb A^1$, then 
we have the following virtual pullback (or Lagrangian pullback) morphism:
\[
\sqrt{g^!}:\CH_*(\mathsf w^{-1}(0))\to\CH_{*+\frac{1}{2}\rank\bb T_{\mathsf{M}/\mathsf{B}}}(\mathsf{M})
\]
with several canonical properties. The isotropic condition means that $g$
is an \textit{exact} symplectic fibration, that is, $\mathsf w=0$; see \cite[Rmk.~5.2.5]{park24}
for details.
We also refer to \cite{park24} and \cite[\S~1]{PY26} for modern introduction about 
the shifted symplectic geometry.

\section{Orientations and \texorpdfstring{$\mathsf{DT}_4$}{DT4}-invariants of a class of toric CY 4-folds}\label{sect on DT toric 4folds}

\subsection{Fourfold vertex formalism for toric Calabi-Yau \texorpdfstring{$4$}{4}-folds}\label{sect on ordinary K-th vertex}
Here we review the theory of fourfold vertices for Hilbert schemes of points on toric Calabi-Yau $4$-folds, following \cite{CK20,CKM22}, which provides a computational framework for the corresponding $\mathsf{DT}_4$-invariants.

Let $X$ be a toric Calabi-Yau $4$-fold. The Hilbert scheme of $n$ points $\Hilb^n(X)$ is a quasi-projective scheme that carries a $(-2)$-shifted symplectic derived enhancement $\dHilb^n(X)\subset\mathsf{dPerf}(X)_0$ as constructed in \cite{PTVV13}; see also \cite[\S~3.1]{CZZ24}. This induces a $3$-term symmetric obstruction theory
\[\bb E_X:=\bb L_{\dHilb^n(X)}|_{\Hilb^n(X)}=\RR\pi_*\RR\scr Hom(\scr I,\scr I)_0[3]\to\mathbf L_{\Hilb^n(X)},\]
where $\pi:X\times\Hilb^n(X)\to \Hilb^n(X)$ is the projection, $\scr I$ is the universal ideal sheaf of the universal substack
$\mcal Z\subset X\times\Hilb^n(X)$ and $(-)_0$ denotes the trace-free part. Fix a choice of orientation\footnote{Here we note that it is unknown whether the general quasi-projective Calabi-Yau 4-folds admit the orientations as initially claimed in \cite{CGJ20,Bojko21}, due to technical issues in the original claims. We refer to the corrected erratum \cite{CGJerratum} of \cite{CGJ20} and related developments in \cite{JU25} for further details.}, then $\bb E_X\to\mathbf L_{\Hilb^n(X)}$ satisfies the isotropy condition by \cite[Prop.~4.3]{OT23}.

Let $\bb T\subset(\bb C^*)^4$ be the Calabi-Yau subtorus defined by the relation $t_1t_2t_3t_4=1$, then the induced $\bb T$-action on $X$ preserves the Calabi-Yau volume form. This induces a canonical action on $\Hilb^n(X)$ such that
$\bb E_X$ is $\bb T$-equivariant and $\bb T$-equivariantly self-dual by relative Serre duality. Consequently, we obtain a $\bb T$-equivariant virtual fundamental class
\[[\Hilb^n(X)]_{\bb T}^{\vir}\in\CH_n^{\bb T}(\Hilb^n(X))\]
as defined in Definition \ref{def of virtual pullbacks}.

We observe that $\Hilb^n(X)$ is not proper in general, which prevents a direct definition of $\mathsf{DT}_4$-invariants.
Nevertheless, we have the following result.
\begin{lem}[{\cite[Lem.~2.1]{CK20}}]\label{lem-fix-locus-Hilb-C4}
The fixed locus $\Hilb^n(X)^{\bb T}=\Hilb^n(X)^{(\bb C^*)^4}$ consists of finitely many isolated reduced points as a scheme.
\end{lem}
This fact, combined with the torus localization formula \cite[Thm.~7.3]{OT23}, allows us to define Donaldson-Thomas invariants through Definition \ref{def-DT4-toric-CY4}.
To facilitate the computation of these invariants, we require an explicit description of the fixed locus $\Hilb^n(X)^{\mathbb{T}}$. This leads us to the following key definition:

\begin{defi}\label{def-solid-partition}
A \textit{solid partition} is a collection of finitely many lattice points
$\pi\subset(\bb Z_{\geq0})^4$ such that if any of $(i+1,j,k,l)$, $(i,j+1,k,l)$, $(i,j,k+1,l)$, $(i,j,k,l+1)$
is in $\pi$, then $(i,j,k,l)\in\pi$. We denote the size of a solid partition by $|\pi|$, i.e. the
number of boxes in $\pi$.
\end{defi}

Since $X$ is a toric Calabi-Yau $4$-fold with fixed locus $X^{\bb T}=X^{(\bb C^*)^4}=\{p_1,...,p_t\}$, each fixed point $p_a$ is contained in a unique maximal $(\bb C^*)^4$-invariant affine open subset $\bb C^4\cong U_{a}\subset X$ such that $\{U_a\}_{a=1}^t$ forms an open cover of $X$. 

\begin{lem}\label{lem-parti-fixed-corr}
Choose coordinates on any $U_a$ such that $(\bb C^*)^4$ acts as the standard way. Then there is a bijection
\[\begin{tikzcd}
	{\left\{\pi=\{\pi^{(a)}\}_{a=1}^{t}\text{ collection of solid partitions}:|\pi|:=|\pi^{(1)}|+\cdots+|\pi^{(t)}|=n\right\}} \\
	{[W_{\pi}]\in\Hilb^n(X)^{\bb T}.}
	\arrow[from=1-1, to=2-1]
	\arrow[shift left=3, from=2-1, to=1-1]
\end{tikzcd}\]
For $[W]\in\Hilb^n(X)^{\bb T}$, it is defined by monomial ideals $\{\scr I_W|_{U_a}=I_a\}$ which induce a collections of solid partitions $\{\pi^{(a)}\}_{a=1}^t$ as
\[\pi^{(a)}:=\{(i,j,k,l):x_1^ix_2^jx_3^kx_4^l\notin I_a\}\subset(\bb Z_{\geq0})^4,\]
see also \cite{MNOP1}.
Moreover, fix $[W]\in\Hilb^n(X)^{\bb T}$ corresponding to $\pi=\{\pi^{(a)}\}_{a=1}^{t}$. Choosing the standard action on $U_a$, we have
\[Z_{\pi^{(a)}}:=\Gamma(\scr O_W|_{U_a})=\sum_{(i,j,k,l)\in\pi^{(a)}}t_1^it_2^jt_3^kt_4^l\in K_0^{\bb T}(\text{pt}).\]
\end{lem}

Next, we need to describe the obstruction theory and express $\DTC_{X,\scr L}(m,q)$ by solid partitions, that is, the vertex formalism.
Let $T^{\vir}:=\bb E^{\vee}_X$ be the virtual tangent bundle. By \cite[Eqn.~(111)]{OT23}, there
exists an induced self-dual obstruction theory on the $\bb T$-fixed locus $\Hilb^n(X)^{\bb T}$, with the virtual
tangent bundle $T^{\vir,\bb T}:=(T^{\vir})^{\bb T\text{-}\fix}$. These are concluded by the following theorem.

\begin{thm}[{\cite[~\S3]{CK18}, \cite[~\S2.4.1]{NP2019} and \cite[~\S1]{CKM22}}]\label{toric-vertex-thm}
Let $[W_{\pi}]\in\Hilb^n(X)^{\bb T}$ be a fixed subscheme corresponding to $\pi=\{\pi^{(a)}\}_{a=1}^{t}$.
By the local-to-global spectral sequence, the virtual tangent space is decomposed as
\[T^{\vir}_{[W_{\pi}]}=-\RR\Hom_X(\scr I_{W_{\pi}},\scr I_{W_{\pi}})_0=-\sum_{a=1}^t\RR\Hom_{U_a}(I_a,I_a)_0\in K_0^{\bb T}(\text{pt}).\]
Define $V^{\DT}_{\pi^{(a)}}:=-\RR\Hom_{U_a}(I_a,I_a)_0\cong Z_{\pi^{(a)}}+\overline{Z_{\pi^{(a)}}}-P_{1234}Z_{\pi^{(a)}}\overline{Z_{\pi^{(a)}}}\in K_0^{\bb T}(\text{pt})$.
Their square roots are chosen by $v^{\DT}_{\pi^{(a)}}:=Z_{\pi^{(a)}}-\overline{P_{123}}Z_{\pi^{(a)}}\overline{Z_{\pi^{(a)}}}\in K_0^{\bb T}(\text{pt})$, which gives $V^{\DT}_{\pi^{(a)}}=v^{\DT}_{\pi^{(a)}}+\overline{v^{\DT}_{\pi^{(a)}}}$. 
Consider the vertices with insertions
\[\widetilde{v}^{\DT,\scr L}_{\pi^{(a)}}:=v^{\DT}_{\pi^{(a)}}-y\cdot\RR\Gamma(U_a,\scr L\otimes\scr O_{W_{\pi}}|_{U_a})^{\vee},\] 
then $\widetilde{v}^{\DT,\scr L}_{\pi^{(a)}}$ does not have $(\bb T\times\bb C_m^*)$-fixed terms. Then we have 
\begin{align}\label{K-th vertex chi toric}
\int_{[\Hilb^n(X)]^{\vir}_{\bb T}}e_{\bb T\times\bb C^*_m}((\scr L^{[n]})^{\vee}\otimes e^m)
=\sum_{\substack{\pi=(\pi^{(a)})_{a=1}^t,|\pi|=n\\ \text{ solid partitions}}}\prod_{a=1}^t(-1)^{\sigma_{\pi^{(a)}}}e_{\bb T\times\bb C_m^*}(-\widetilde{v}^{\DT,\scr L}_{\pi^{(a)}}),
\end{align}
where $\sigma_{\pi^{(a)}}\in\bb Z/2$ be the sign determined by the orientation.
We also denote $\widetilde{v}^{\DT}_{\pi^{(a)}}:=\widetilde{v}^{\DT,\scr O}_{\pi^{(a)}}$.
Here, for $f(t_1,...,t_4)\in K_0^{\bb T}(\text{pt})$, we define $\overline{f(t_1,...,t_4)}=f(t_1^{-1},...,t_4^{-1})$
and $P_{123}=(1-t_1)(1-t_2)(1-t_3)$ and $P_{1234}$ similarly.
\end{thm}

\subsection{Canonical orientation of Hilbert schemes of a class of toric CY \texorpdfstring{$4$}{4}-folds}
In this section, we will consider the orientations, sign rules, and Donaldson-Thomas invariants of $X=Y\times\bb C$ for some toric Calabi-Yau $3$-fold $Y$.

We first construct a canonical orientation in a more general situation for subsequent applications as follows.

\begin{lem}[{\cite[Thm.~4.3,~Prop.~B.2]{CZZ24}} and {\cite[~\S 4]{KR}}]\label{lem-ori-identification}
Let $\mcal Y$ be a quasi-projective smooth Deligne-Mumford stack of dimension $3$. Identification $\mcal X:=\Tot(\omega_{\mcal Y})$ induces a canonical orientation of the tangent complex $\bb T_{\dHilb^{\alpha}(\mcal X)}$
and, hence, an orientation of the symmetric obstruction theory of $\Hilb^{\alpha}(\mcal X)$ for any $\alpha\in F_0K_{\text{cpt}}^{\text{num}}(\mcal X)$; see \cite[Def.~2.13]{CZZ24} for precise definitions.

Moreover, if $\mcal X=\bb C^4, \mcal Y=\bb C^3$ and $\alpha=n[\scr O_{\text{pt}}]$, this orientation determined the sign rule
$\sigma_{\pi}\equiv\mu_{\pi}\mmod{2}$ in Eqn. (\ref{K-th vertex chi toric}) where $\mu_{\pi}:=|\{(i,i,i,j)\in\pi:j>i\}|$.
\end{lem}
\begin{proof}
The orientation induced by identification $\mcal X=\Tot(\omega_{\mcal Y})$ is given in \cite[Thm.~4.3]{CZZ24}. The sign rules are computed in \cite[~\S 4]{KR} and compared with this orientation in \cite[Prop.~B.2]{CZZ24}; see also the proofs in $\S$\ref{section-compare-oris}.
Here we recall the construction of the orientation discussed in \cite[Thm.~4.3]{CZZ24} and we will use it later. 

For $\bullet=\mcal X$ or $\mcal Y$, let $\textsf{M}_{\alpha}(\bullet)$ be the moduli stack of zero-dimensional sheaves for $\alpha\in F_0K_{\text{cpt}}^{\text{num}}(\bullet)$ and let $\bb F_{\bullet}\to\textsf{M}_n(\bullet)\times\bullet$ be the universal object
with projection $\pi_{M_{\bullet}}:\textsf{M}_{\alpha}(\bullet)\times\bullet\to\textsf{M}_{\alpha}(\bullet)$.
The canonical projection $\pi:\mcal X\to\mcal Y$ induces $\pi_{\dagger}:\mathsf M_{\alpha}(\mcal X)\to\mathsf M_{\alpha}(\mcal Y)$ and its derived enhancement, which we will use the same notation $\pi_{\dagger}:\dMSh_{\alpha}(\mcal X)\to\dMSh_{\alpha}(\mcal Y)$.

Consider the forgetful map $f:\Hilb^{\alpha}(\mcal X)\to\textsf{M}_{\alpha}(\mcal X)$ given by $(\scr O\twoheadrightarrow\scr F)\mapsto\scr F$.
By some homological algebra (see \cite[~\S 3.2,~Thm.~6.3]{CMT21}), we have 
\begin{equation*}
\begin{tikzcd}[row sep=small]
	{\det(\bb T_{\dHilb^{\alpha}(\mcal X)}|_{\Hilb^{\alpha}(\mcal X)})} & {f^*\det(\bb T_{\dMSh_{\alpha}(\mcal X)}|_{\MSh_{\alpha}(\mcal X)})} \\
	{\det(\RR\pi_*\RR\scr Hom(\scr I,\scr I)_0[1])} & {f^*\det(\RR\pi_{M_{\mcal X},*}\RR\scr Hom(\bb F_{\mcal X}, \bb F_{\mcal X})[1]).}
	\arrow["\cong", from=1-1, to=1-2]
	\arrow[no head, from=1-1, to=2-1]
	\arrow[shift left, no head, from=1-1, to=2-1]
	\arrow[no head, from=1-2, to=2-2]
	\arrow[shift left, no head, from=1-2, to=2-2]
	\arrow["\cong", from=2-1, to=2-2]
\end{tikzcd}
\end{equation*}
This shows that an orientation of $\textsf{dM}_{\alpha}(\mcal X)$ will induce an orientation of $\dHilb^{\alpha}(\mcal X)$. In what follows, we need to find an orientation of $\textsf{dM}_{\alpha}(\mcal X)$.

The first approach follows from $\dMSh(\mcal X)\cong T^*[-2]\dMSh(\mcal Y)$. However, we will use the second approach.
By \cite[Eqn.~(4.8),~(4.17)]{CZZ24} and some base change, we have the following distinguished triangle:
\begin{footnotesize}
    \begin{equation}\label{eqn-isotropic-quotient-orientation}
    \pi_{\dagger}^*(\RR\pi_{M_{\mcal Y},*}\RR\scr Hom(\bb F_{\mcal Y}, \bb F_{\mcal Y}))^{\vee}[-3]\to
    \RR\pi_{M_{\mcal X},*}\RR\scr Hom(\bb F_{\mcal X}, \bb F_{\mcal X})[1]
    \to\pi_{\dagger}^*\RR\pi_{M_{\mcal Y},*}\RR\scr Hom(\bb F_{\mcal Y}, \bb F_{\mcal Y})[1]
\end{equation}
\end{footnotesize}
of tangent complexes induced by $\pi_{\dagger}:\dMSh_{\alpha}(\mcal X)\to\dMSh_{\alpha}(\mcal Y)$.
This gives us an orientation
\begin{equation*}
   \det(\RR\pi_{M_{\mcal X},*}\RR\scr Hom(\bb F_{\mcal X}, \bb F_{\mcal X})[1])
    \cong\frac{\det(\pi_{\dagger}^*\RR\pi_{M_{\mcal Y},*}\RR\scr Hom(\bb F_{\mcal Y}, \bb F_{\mcal Y})[1])}{\det(\pi_{\dagger}^*\RR\pi_{M_{\mcal Y},*}\RR\scr Hom(\bb F_{\mcal Y}, \bb F_{\mcal Y})[1])}
    \cong\scr O.
\end{equation*}
See \cite[Rmk.~4.4]{CZZ24} for the precise choice of orientations.
\end{proof}

Here we give the sign rules determined by this orientation via the compatibility of the orientations of toric charts.

\begin{prp}\label{thm-sign-rule-toric}
Let $X=Y\times\bb C$ where $Y$ is a toric Calabi-Yau $3$-fold, the orientation defined in Lemma \ref{lem-ori-identification} determines
\begin{align}\label{sign rule toric}
\int_{[\Hilb^n(X)]^{\vir}_{\bb T}}e_{\bb T\times\bb C^*_m}((\scr L^{[n]})^{\vee}\otimes e^m)
=\sum_{{\substack{\pi=(\pi^{(a)})_{a=1}^t,|\pi|=n\\ \text{ solid partitions}}}}\prod_{a=1}^t(-1)^{\mu_{\pi^{(a)}}}e_{\bb T\times\bb C_m^*}(-\widetilde{v}^{\DT,\scr L}_{\pi^{(a)}}),
\end{align}
where $\mu_{\pi^{(a)}}=|\{(i,i,i,j)\in\pi^{(a)}:j>i\}|$ which coincides with the sign rules in  \cite[Rmk.~1.18]{CKM22}.
\end{prp}
\begin{proof}
As before, write $Y^{(\bb C^*)^3}=\{q_1,...,q_t\}$ and let $V_a\cong\bb C^3$ be the maximal affine toric chart containing $q_a$.
Let $\pi:X=\Tot(\omega_Y)\to Y$ be the projection and set $U_a:=\pi^{-1}(V_a)=\Tot(\omega_{V_a})$, with projection $\pi_a:U_a\to V_a$.
Choose toric coordinates on $V_a$ and the induced canonical bundle coordinate $x_4$ on $U_a$.
Thus these local presentations are compatible with $X=\Tot(\omega_Y)$ and its Calabi-Yau form, and $X^{\bb T}=\{p_1,...,p_t\}$ where $p_a=q_a\times\{0\}$.

For a fixed point $[Z]\in\Hilb^n(X)^{\bb T}$, write $Z=\coprod_aZ_a$, where $Z_a\subset U_a$ is supported at $p_a$ and corresponds to the solid partition $\pi^{(a)}$. 
Put $n_a:=\operatorname{length}(Z_a)$, so $\sum_an_a=n$.
Put
\[
\scr F:=\scr O_Z,\quad \scr F_a:=\scr O_{Z_a},\quad \scr G:=\pi_*\scr F,\quad \scr G_a:=(\pi_a)_*\scr F_a.
\]
Here empty $Z_a$ are allowed, and we also regard $\scr F_a$ and $\scr G_a$ as sheaves on $X$ and $Y$, respectively, by extension by zero.
Then $\scr F=\bigoplus_a\scr F_a$ and $\scr G=\bigoplus_a\scr G_a$.

By Lemma \ref{lem-ori-identification}, we first compare the orientations on the moduli stacks of zero-dimensional sheaves. Denote
\[
\begin{aligned}
\bb D&:=\RR\Hom_X(\scr F,\scr F)[1],&\quad \bb V_1&:=\RR\Hom_Y(\scr G,\scr G)[1],\\
\bb D_a&:=\RR\Hom_{U_a}(\scr F_a,\scr F_a)[1],&\quad \bb V_{1,a}&:=\RR\Hom_{V_a}(\scr G_a,\scr G_a)[1].
\end{aligned}
\]
Since the points $p_a$ and their projections $q_a$ are respectively distinct, we have
\[
\RR\Hom_X(\scr F_a,\scr F_b)=0,\quad \RR\Hom_Y(\scr G_a,\scr G_b)=0\quad\text{for }a\ne b.
\]
Hence restriction gives $\bb T$-equivariant quasi-isomorphisms
\[
\bb D\cong\bigoplus_a\bb D_a,\quad \bb V_1\cong\bigoplus_a\bb V_{1,a}.
\]
The spectral construction in Eqn. (\ref{eqn-isotropic-quotient-orientation}) is compatible with restriction to these charts and with the above decompositions.
Moreover, the Grothendieck-Serre trace is the sum of the traces on the disjoint supports. Thus we obtain the following commutative diagram of distinguished triangles and isotropic quotients:
\begin{equation}\label{diagram-ori-torus-3}
\begin{tikzcd}
\bb V_1^{\vee}[-2] \arrow[r,"\cong"] \arrow[d] & \bigoplus_a\bb V_{1,a}^{\vee}[-2] \arrow[d] \\
\bb D \arrow[r,"\cong"] \arrow[d] & \bigoplus_a\bb D_a \arrow[d] \\
\bb V_1 \arrow[r,"\cong"] & \bigoplus_a\bb V_{1,a}.
\end{tikzcd}
\end{equation}
Here the top horizontal arrow is the inverse dual of the bottom one, shifted by $[-2]$, and the middle arrow preserves the Serre pairings.

For determinant evaluations we use the duality conventions in \cite[~\S 4.1]{OT23}; orientations and their transfer are understood in the conventions of \cite[Eqn.~(4.4),~Rmk.~4.4]{CZZ24}.
Let
\[
p_{\bb V_1}:\det\bb D\cong\det\bb V_1\otimes\det\bb V_1^{\vee}\cong\bb C
\]
be the trivialization induced by the isotropic quotient, and define $p_{\bb V_{1,a}}$ similarly.
By Serre duality on the Calabi-Yau $3$-folds $V_a$, we have
$\rank\bb V_{1,a}=0=\rank\bb V_1$.
Therefore, the orientations specified in \cite[Rmk.~4.4]{CZZ24} are
\begin{equation}\label{eqn-orientation-1-2}
o_1:=(-\sqrt{-1})^{\rank\bb V_1}p_{\bb V_1}^{-1}=p_{\bb V_1}^{-1},\quad
o_{1,a}:=p_{\bb V_{1,a}}^{-1}.
\end{equation}
By Figure (\ref{diagram-ori-torus-3}) and compatibility of determinants with direct sums, these satisfy
\[
o_1=\bigotimes_a o_{1,a}
\]
under $\det\bb D\cong\bigotimes_a\det\bb D_a$.

Now we compare the induced orientations on the Hilbert schemes. Let $\scr I_Z\subset\scr O_X$ and $\scr I_a\subset\scr O_{U_a}$ be the ideal sheaves of $Z$ and $Z_a$, and denote
\[
\bb E:=\RR\Hom_X(\scr I_Z,\scr I_Z)_0[1],\quad \bb E_a:=\RR\Hom_{U_a}(\scr I_a,\scr I_a)_0[1].
\]
Here $\RR\scr Hom_X(\scr I_Z,\scr I_Z)_0$ denotes the homotopy fibre of the trace map to $\scr O_X$.
Its cohomology sheaves are supported on $Z$, since $\scr I_Z=\scr O_X$ away from $Z$ and the trace map there is an isomorphism.
By the local-to-global spectral sequence, restriction induces a quasi-isomorphism
\begin{equation}\label{Eqn-tangentcomplex-many-Hilbs}
\Res_{[Z]}:\bb E\xrightarrow{\cong}\bigoplus_a\bb E_a,
\end{equation}
which preserves the Serre pairings by compatibility of the compactly supported trace with restriction.

As in \cite[Eqn.~(4.4)]{CZZ24}, the forgetful maps from Hilbert schemes to moduli stacks of sheaves give determinant identifications
\[
\eta:\det\bb E\xrightarrow{\cong}\det\bb D,\quad
\eta_a:\det\bb E_a\xrightarrow{\cong}\det\bb D_a,
\]
under which orientations are transferred as in Lemma \ref{lem-ori-identification}.
We check their compatibility with the above decompositions. The identification $\eta$ is constructed from the distinguished triangles
\[
\begin{aligned}
\RR\Hom_X(\scr I_Z,\scr F)&\longrightarrow\bb E\longrightarrow\RR\Hom_X(\scr F,\scr O_X)[2],\\
\RR\Gamma(X,\scr F)&\longrightarrow\RR\Hom_X(\scr I_Z,\scr F)\longrightarrow\bb D
\end{aligned}
\]
and Grothendieck-Serre duality. All terms and maps in these triangles decompose according to the supports of the $\scr F_a$, giving the corresponding triangles on $U_a$.
In particular, put
\[
A:=\RR\Gamma(X,\scr F)=\bigoplus_aA_a,\quad A_a:=\RR\Gamma(U_a,\scr F_a).
\]
Since $\scr F_a$ is zero-dimensional, $A_a$ is concentrated in degree zero and has dimension $n_a$.
The framing terms are $A$ and $\RR\Hom_X(\scr F,\scr O_X)[2]\cong A^{\vee}[-2]$, whose determinant factors are cancelled by the evaluation pairing.
With the above duality conventions, the dual factors are taken in reverse order. For two summands, regrouping the paired determinant factors from $\mathbf a\wedge\mathbf b\wedge\mathbf b^{\vee}\wedge\mathbf a^{\vee}$ to $\mathbf a\wedge\mathbf a^{\vee}\wedge\mathbf b\wedge\mathbf b^{\vee}$ gives the sign $(-1)^{2n_an_b}=1$.
The remaining determinant factors $\det\bb D_a$ have degree zero, since $\rank\bb D_a=0$.
Consequently, the determinant identifications, including the transfer of orientations, are compatible with direct sums.

Thus we have the following commutative diagram:
\begin{equation}\label{diagram-ori-torus-5}
\begin{tikzcd}
\bb C \arrow[r,"o_1"] \arrow[d,equal] & \det\bb D \arrow[r,"\eta^{-1}"] \arrow[d,"\cong"] & \det\bb E \arrow[d,"\det\Res_{[Z]}"] \\
\bb C \arrow[r,"\bigotimes_a o_{1,a}"'] & \bigotimes_a\det\bb D_a \arrow[r,"\bigotimes_a\eta_a^{-1}"'] & \bigotimes_a\det\bb E_a.
\end{tikzcd}
\end{equation}
In particular, the canonical orientation of $\bb E$ agrees with the product of the canonical orientations of $\bb E_a$.

Finally, by Eqn. (\ref{Eqn-tangentcomplex-many-Hilbs}), the virtual tangent character and its chosen square root in Theorem \ref{toric-vertex-thm} are
\[
T^{\vir}_{\Hilb^n(X),[Z]}=\sum_aV^{\DT}_{\pi^{(a)}},\quad
\sqrt{T^{\vir}_{\Hilb^n(X),[Z]}}=\sum_av^{\DT}_{\pi^{(a)}}
\]
which coincides the canonical one in Theorem \ref{toric-vertex-thm}.
For each $U_a=\Tot(\omega_{V_a})\cong\bb C^4$, the second statement of Lemma \ref{lem-ori-identification} gives the sign $(-1)^{\mu_{\pi^{(a)}}}$ with preferred axis $x_4$.
The tautological insertion also decomposes, since
\[
\RR\Gamma(X,\scr L\otimes \scr F)\cong\bigoplus_a\RR\Gamma(U_a,\scr L|_{U_a}\otimes \scr F_a).
\]
Hence the contribution of $[Z]$ is
\[
\prod_a(-1)^{\mu_{\pi^{(a)}}}e_{\bb T\times\bb C_m^*}(-\widetilde v^{\DT,\scr L}_{\pi^{(a)}}).
\]
Summing over the fixed points in Theorem \ref{toric-vertex-thm}, we obtain Eqn. (\ref{sign rule toric}).
\end{proof}

\begin{nota}
The compatibility of the sign rules established in Proposition \ref{thm-sign-rule-toric} also complements the proof of \cite[Thm.~5.17]{CKM23}.
\end{nota}
 
\begin{nota}
We can take other square roots by choosing other preferred $x_i$ axes, that is, define
\[
v_{\pi^{(a)}}^{\DT,i}:=Z_{\pi^{(a)}}-\overline{P_{jkl}}Z_{\pi^{(a)}}\overline{Z_{\pi^{(a)}}},
\]
where $\{i,j,k,l\}=\{1,2,3,4\}$. In \cite[Thm.~2.8]{monavari22}, the author showed that the sign rules are canonical, that is, one has
$(-1)^{\mu_{\pi^{(a)}}^i}e_{\bb T}(-\widetilde{v}^{\DT,i}_{Z_a})=(-1)^{\mu_{\pi^{(a)}}^j}e_{\bb T}(-\widetilde{v}^{\DT,j}_{Z_a})$.
Therefore,
\begin{align*}
\int_{[\Hilb^n(X)]^{\vir}_{\bb T}}e_{\bb T\times\bb C^*_m}((\scr L^{[n]})^{\vee}\otimes e^m)
=\sum_{{\substack{\pi=(\pi^{(a)})_{a=1}^t,|\pi|=n\\ \text{ solid partitions}}}}\prod_{a=1}^t(-1)^{\mu_{\pi^{(a)}}^i}e_{\bb T\times\bb C_m^*}(-\widetilde{v}^{\DT,i,\scr L}_{\pi^{(a)}}),
\end{align*}
where $\widetilde{v}^{\DT,i,\scr L}_{\pi^{(a)}}:=v^{\DT,i}_{\pi^{(a)}}-y\cdot\RR\Gamma(U_a,\scr L\otimes\scr O_{Z_{\pi}}|_{U_a})^{\vee}$ and
$\mu_{\pi^{(a)}}^i=|\{(a_1,a_2,a_3,a_4)\in\pi^{(a)}:a_j=a_k=a_l<a_i\}|$.
For example, let $X=Z\times\bb C^2$ where $Z$ is a toric Calabi-Yau surface, then two choices of identifications
$X=\Tot(\omega_{Z\times\bb C_3})$ and $\Tot(\omega_{Z\times\bb C_4})$ induce the same orientations.
\end{nota}

\subsection{Zero-dimensional \texorpdfstring{$\mathsf{DT}_4$}{DT4}-invariants of a class of toric CY \texorpdfstring{$4$}{4}-folds}
As an application, we can deduce $\mathsf{DT}_4$-invariants
and confirm Conjecture \ref{conj-torus-DT4-general} for the following toric CY $4$-folds.
\begin{thm}\label{coro-DT-specil-toric}
Let $X=Y\times\bb C$ where $Y$ is a toric Calabi-Yau $3$-fold. Fix the orientation determined by Lemma \ref{lem-ori-identification}, then we have
\begin{equation}
    \DTC_{X,\scr L}(m,q)=M(-q)^{\int_Xc_1^{\bb T}(\scr L\otimes y^{-1})c_3^{\bb T}(X)}.
\end{equation}
\end{thm}
\begin{proof}
Applying \cite[Thm.~6.8]{CZZ24} to $\scr L\otimes y^{-1}$, we obtain
\begin{equation}
    \DTC_{\bb C^4,\scr L}(m,q)=M(-q)^{\int_{\bb C^4}c_1^{\bb T}(\scr L\otimes y^{-1})c_3^{\bb T}(\bb C^4)}.
\end{equation}
Combining this with the vertex formula and the sign rules in Eqn. (\ref{sign rule toric}) of Proposition \ref{thm-sign-rule-toric}, we have
\begin{align*}
    \DTC_{X,\scr L}(m,q)&=\sum_{\substack{\pi=(\pi^{(a)})_{a=1}^t\\ \text{ solid partitions}}}q^{|\pi|}\prod_{a=1}^t(-1)^{\mu_{\pi^{(a)}}}e_{\bb T\times\bb C_m^*}(-\widetilde{v}^{\DT,\scr L}_{\pi_a})\\
    &=\prod_{a=1}^t\sum_{\substack{\pi^{(a)}\\\text{ solid partition}}}q^{|\pi^{(a)}|}(-1)^{\mu_{\pi^{(a)}}}e_{\bb T\times\bb C_m^*}(-\widetilde{v}^{\DT,\scr L}_{\pi_a})\\
    &=\prod_{a=1}^t\DTC_{U_a,\scr L|_{U_a}}(m,q)=M(-q)^{\sum_{a=1}^t\int_{U_a}c_1^{\bb T}(\scr L|_{U_a}\otimes y^{-1})c_3^{\bb T}(U_a)}.
\end{align*}
Using Atiyah-Bott localization formula twice and $T_{p_a}X=T_{p_a}U_a$, one has
\begin{align*}
    &\int_Xc_1^{\bb T}(\scr L\otimes y^{-1})c_3^{\bb T}(X)
    =\sum_{a=1}^t\int_{p_a}\frac{\iota_a^*(c_1^{\bb T}(\scr L\otimes y^{-1})c_3^{\bb T}(X))}{e^{\bb T}(T_{p_a}X)}
    =\sum_{a=1}^t\int_{U_a}c_1^{\bb T}(\scr L|_{U_a}\otimes y^{-1})c_3^{\bb T}(U_a),
\end{align*}
where $\iota_a:\{p_a\}\hookrightarrow X$ be the inclusion.
Then the result follows.
\end{proof}
As an immediate consequence of Theorem \ref{coro-DT-specil-toric}, we obtain the following special case that plays a crucial role in computing the Donaldson-Thomas invariants of $[\bb C^4/\bb Z_r]$.

\begin{coro}\label{coro-DT-ArC2}
Let $X=\widetilde{\bb C^2/\bb Z_r}\times\bb C^2$
where $\widetilde{\bb C^2/\bb Z_r}$ be the minimal crepant resolution of $A_{r-1}$-singularity 
$\bb C^2/_{\text{aff}}\bb Z_r$ with action $\zeta_r\cdot(x_1,x_2)=(\zeta_r x_1,\zeta_r^{-1}x_2)$ for $\zeta_r=e^{2\pi\sqrt{-1}/r}$. We fix the orientation determined by 
Lemma \ref{lem-ori-identification}, then 
\begin{equation}
    \DTC_{\widetilde{\bb C^2/\bb Z_r}\times\bb C^2,\scr O}(m,q)=M(-q)^{-\frac{m}{s_4}
    \left(\frac{r(s_1+s_2)}{s_3}+\frac{(s_1+s_2)(s_1+s_2+s_3)}{rs_1s_2}\right)}.
\end{equation}
\end{coro}

\section{Orientations and sign rules of the Hilbert scheme \texorpdfstring{$\Hilb^{R}([\mathbb C^4/\bb Z_r])$}{Hilb R [C4/Zr]}}\label{sect on ori of C4/Zr}

\subsection{Fourfold vertex formalism for \texorpdfstring{$\Hilb^{R}([\mathbb C^4/\bb Z_r])$}{Hilb R [C4/Zr]}}
In this section, we adapt the general framework developed in \cite{CKM23} to our specific case, the Hilbert scheme $\Hilb^{R}([\mathbb C^4/\bb Z_r])$.

For a finite subgroup $G<\text{SU}(4)$ acting on $\bb C^4$,
the K-theoretic classes in $F_0K_{\text{cpt}}^{\text{num}}([\bb C^4/G])$ are equivalent to the finite-dimensional representations of $G$.
For any finite-dimensional $G$-representation $R$, the stack $\Hilb^{R}([\bb C^4/G])$ is a quasi-projective scheme (this holds for more general Deligne-Mumford stacks; see \cite{OS03}) given by
\begin{displaymath}
	\Hilb^{R}([\bb C^4/G])=\left\{[W]\in\Hilb^{\dim R}(\bb C^4)\left|
		\begin{aligned}
		&W\subset\bb C^4\text{ is 0-dim. }G\text{-invariant}\\
		&\text{with }H^0(W,\scr{O}_W)\cong R\text{ as }G\text{-reps}.
		\end{aligned}\right.
		\right\}.
\end{displaymath}
This gives a canonical decomposition:
\[\Hilb^n(\bb C^4)^{G}=\coprod_{\dim R=n}\Hilb^{R}([\bb C^4/G]).\]

Return to the case $G=\bb Z_r<\text{SU}(2)<\text{SU}(4)$ with the action on $\bb C^4$ given by 
\[\zeta_r\cdot(x_1,x_2,x_3,x_4)=(\zeta_r x_1,\zeta^{-1}_rx_2,x_3,x_4).\]
We fix a finite-dimensional representation $R=\rho_0^{n_0}\oplus\cdots\oplus
\rho_{r-1}^{n_{r-1}}$ of $\bb Z_r$, where $\rho_k$ is the $1$-dimensional representation induced by $\zeta_r^k$. 
The Hilbert scheme $\Hilb^{R}([\bb C^4/\bb Z_r])$ is a quasi-projective scheme with a natural $(-2)$-shifted 
symplectic derived enhancement
$\dHilb^{R}([\bb C^4/\bb Z_r])$, which induces a similar $\bb T$-equivariant $3$-term symmetric obstruction 
theory
\begin{equation}\label{eqn-obs-theory-C4Zr}
    \bb E:=\bb L_{\dHilb^{R}([\bb C^4/\bb Z_r])}|_{\Hilb^{R}([\bb C^4/\bb Z_r])}=\RR\pi_*\RR\scr Hom(\scr I,\scr I)_0[3]\to\mathbf L_{\Hilb^R([\bb C^4/\bb Z_r])}.
\end{equation}
Fix an orientation
$o$, which will be constructed later, then it satisfies the isotropy condition by \cite[Prop.~4.3]{OT23}.
Therefore, we obtain the $\bb T$-equivariant virtual fundamental class
\[[\Hilb^{R}([\bb C^4/\bb Z_r])]_{\bb T,o}^{\vir}\in\CH_{n_0}^{\bb T}(\Hilb^{R}([\bb C^4/\bb Z_r]))\]
by Definition \ref{def of virtual pullbacks} again.

As before, the fixed locus $\Hilb^R([\bb C^4/\bb Z_r])^{\bb T}$ can also be described by solid partitions; see \cite[~\S 5.6]{CKM23} for the description for more general $G$.
\begin{defi}\label{def-solid-partition-coloured}
Fix a finite-dimensional representation $R=\rho_0^{n_0}\oplus\cdots\oplus
\rho_{r-1}^{n_{r-1}}$ of $\bb Z_r$.
\begin{enumerate}
    \item For a solid partition $\pi$, a box $(i,j,k,l)\in \pi$ is called \textit{$\rho_s$-colored} if $i-j\equiv s\mmod r$. We write $|\pi|_{\rho_s}$ as the number of boxes in $\pi$ of color $\rho_s$.
    \item A solid partition $\pi$ is called \textit{$R$-colored} if
    $|\pi|_{\rho_s}=n_s$ for any $s=0,...,r-1$. In this case, we write $|\pi|_{\bb Z_r}=R$.
\end{enumerate}
\end{defi}

\begin{lem}
By a similar construction as Lemma \ref{lem-parti-fixed-corr}, $\Hilb^R([\bb C^4/\bb Z_r])^{\bb T}$ consists of finite isolated reduced points as a scheme and we have a
bijection
\[\Hilb^R([\bb C^4/\bb Z_r])^{\bb T}\leftrightarrow\{\pi\text{ solid partition}:|\pi|_{\bb Z_r}=R\}.\]
\end{lem}

Similarly, we can define the invariants as before.

\begin{defi}[$\DT_4$-invariants of {$[\bb C^4/\bb Z_r]$}]\label{def-DT4-invariants-C4Zr}
Let $\scr L$ be a $\bb T$-equivariant line bundle on $[\bb C^4/\bb Z_r]$.
The \textit{tautological bundle of $\scr L$} on $\Hilb^R([\bb C^4/\bb Z_r])$ given by
\[\scr L^{[R]}:=\RR\pi_*(\pi_{[\bb C^4/\bb Z_r]}^*\scr L\otimes\scr O_{\mcal Z}),\]
where $\mcal Z\subset [\bb C^4/\bb Z_r]\times\Hilb^R([\bb C^4/\bb Z_r])$ is the universal closed subscheme and $\pi_{[\bb C^4/\bb Z_r]}, \pi$
are projections to the first and second factors of $[\bb C^4/\bb Z_r]\times\Hilb^R([\bb C^4/\bb Z_r])$, respectively. Let $\bb C^*_m$ be an $1$-dimensional torus with trivial action on $[\bb C^4/\bb Z_r]$ and let $y=e^m$ be a character of $\bb C_m^*$
with weight $m=c_1(y)$. We define
the generating function of \textit{$\DT_4$-invariants with tautological insertions} 
as
    \begin{align*}
    &\DTC_{[\bb C^4/\bb Z_r],\scr L,o}(m,q_0,...,q_{r-1})\\
    &:=1+\sum_{n\geq1}\sum_{n_0+\cdots+ n_{r-1}=n\leftrightarrow R}
    \left(\int_{[\Hilb^{R}([\bb C^4/\bb Z_r])]^{\vir}_{\bb T,o}}e_{\bb T\times\bb C^*_m}((\scr L^{[R]})^{\vee}\otimes e^m)\right)\cdot q_0^{n_0}\cdots q_{r-1}^{n_{r-1}}\\
    &\in\text{Frac}\CH^*_{\bb T}(\text{pt})\otimes\CH^*_{\bb C^*_m}(\text{pt})[[q_0,...,q_{r-1}]].
\end{align*}
Here $n_0+\cdots+ n_{r-1}=n\leftrightarrow R$ means that
$R=\rho_0^{n_0}\oplus\cdots\oplus\rho_{r-1}^{n_{r-1}}$ with $\dim R=n$.
\end{defi}

Next, we need to use the vertex formalism to describe the obstruction theory defined above again. 
Denoting $T^{\vir}_{[\bb C^4/\bb Z_r]}:=\bb E^{\vee}$ the virtual tangent bundle, we have $T^{\vir}_{[\bb C^4/\bb Z_r]}=(T^{\vir}_{\bb C^4})^{\bb Z_r\text{-}\fix}$. Similarly to the toric case, there
exists an induced self-dual obstruction theory of $\Hilb^R([\bb C^4/\bb Z_r])^{\bb T}$ with the virtual
tangent bundle $T^{\vir,\bb T}_{[\bb C^4/\bb Z_r]}:=(T^{\vir}_{[\bb C^4/\bb Z_r]})^{\bb T\text{-}\fix}$.

Similarly to $\S$\ref{sect on ordinary K-th vertex},
for $\pi\leftrightarrow[W]\in\Hilb^R([\bb C^4/\bb Z_r])^{\bb T}$ one has
\[Z_{\pi}=H^0(W,\scr O_W)=\sum_{(i,j,k,l)\in\pi}\rho_{i-j}t_1^it_2^jt_3^kt_4^l
=\sum_{\ell=0}^{r-1}\rho_{\ell}Z_{\pi}^{(\ell)}
\in K_0^{\bb T\times\bb Z_r}(\text{pt})\]
as $\bb T\times\bb Z_r$-modules. Here $P_{123}=(1-\rho_1t_1)(1-\rho_{-1}t_2)(1-t_3)$ and $P_{1234}=P_{123}(1-t_4)$, with duals also inverting the $\bb Z_r$-characters.
Then $T^{\vir}_{\bb C^4,[Z]}=Z_{\pi}+\overline{Z_{\pi}}-P_{1234}Z_{\pi}\overline{Z_{\pi}}$.
Hence, 
\[T^{\vir}_{[\bb C^4/\bb Z_r],[W]}=(T^{\vir}_{\bb C^4,[W]})^{\bb Z_r\text{-}\fix}
=v^{\DT,\bb Z_r}_{\pi}+\overline{v^{\DT,\bb Z_r}_{\pi}},\] with the choice of the canonical square roots\footnote{Note that these calculations of vertices are also appeared as an example in \cite[Exam.~5.10]{CKM23}. However, they have some small typos. It should be $l-k\equiv\pm1\mmod r$ instead of $l+k\equiv\pm1\mmod r$ as originally stated.}
\begin{equation}
    \begin{aligned}
    v^{\DT,\bb Z_r}_{\pi}&=(v^{\DT}_{\pi})^{\bb Z_r\text{-}\fix}
=\left(Z_{\pi}-\overline{P_{123}}Z_{\pi}\overline{Z_{\pi}}\right)^{\bb Z_r\text{-}\fix}\\
&= Z_{\pi}^{(0)}-(1-t_3^{-1})\left((1+t_1^{-1}t_2^{-1})\sum_{l\equiv k\mmod r}Z_{\pi}^{(l)}\overline{Z_{\pi}^{(k)}}\right.\\
&\left.-t_1^{-1}\sum_{l-k\equiv 1\mmod r}Z_{\pi}^{(l)}\overline{Z_{\pi}^{(k)}}-t_2^{-1}\sum_{l-k\equiv -1\mmod r}Z_{\pi}^{(l)}\overline{Z_{\pi}^{(k)}}
\right).
    \end{aligned}
\end{equation}
Vertices with insertions are defined as
\begin{equation}
    \begin{aligned}
        \widetilde{v}^{\DT,\bb Z_r}_{\pi}&=(\widetilde v^{\DT}_{\pi})^{\bb Z_r\text{-}\fix}
= Z_{\pi}^{(0)}-\overline{Z_{\pi}^{(0)}}y-(1-t_3^{-1})\left((1+t_1^{-1}t_2^{-1})\sum_{l\equiv k\mmod r}Z_{\pi}^{(l)}\overline{Z_{\pi}^{(k)}}\right.\\
&\left.-t_1^{-1}\sum_{l-k\equiv 1\mmod r}Z_{\pi}^{(l)}\overline{Z_{\pi}^{(k)}}-t_2^{-1}\sum_{l-k\equiv -1\mmod r}Z_{\pi}^{(l)}\overline{Z_{\pi}^{(k)}}
\right).
    \end{aligned}
\end{equation}
Parallel to Theorem \ref{toric-vertex-thm}, we have
\begin{equation}
\DTC_{[\bb C^4/\bb Z_r],\scr O,o}(y,q_0,...,q_{r-1})
=\sum_{\pi}(-1)^{\sigma_{\pi}^{\bb Z_r}}e_{\bb T\times\bb C_m^*}(-\widetilde{v}^{\DT,\bb Z_r}_{\pi})
q_0^{|\pi|_{\rho_0}}\cdots q_{r-1}^{|\pi|_{\rho_{r-1}}}
\end{equation}
for some sign rules $\sigma_{\pi}^{\bb Z_r}\in\bb Z/2$ determined by the choice of orientation $o$.
The main aim of this section is to find the sign rules $\sigma_{\pi}^{\bb Z_r}$ for some canonical choice of orientations.

\subsection{Non-commutative moduli spaces over \texorpdfstring{$[\mathbb C^4/\bb Z_r]$}{[C4/Zr]} via quiver representations}
In this section, we construct the non-commutative Hilbert scheme
$\ncHilb^{R}([\bb C^4/\bb Z_r])$ and the non-commutative 
moduli stack of zero dimensional sheaves $\ncM_{R}([\bb C^4/\bb Z_r])$.
These spaces will play a crucial role in our computations of sign rules.

\subsubsection{Non-commutative Hilbert scheme \texorpdfstring{$\ncHilb^{R}([\bb C^4/\bb Z_r])$}{ncHilb R [C4/Zr]}}\label{sec-noncommu-Hilb-cons}

Consider the following framed quiver

\begin{figure}[!ht]
\begin{tikzcd}
	&& {W_0\diamond} \\
	{\widetilde{Q}:=} &&& {\bullet V_1} & \cdots & {\bullet V_{r-2}} \\
	&& {\bullet V_0} &&&& {\bullet V_{r-1},}
	\arrow["u"', from=1-3, to=3-3]
	\arrow["{C_{11}}", from=2-4, to=2-4, loop, in=60, out=120, distance=5mm]
	\arrow["{C_{10}}", from=2-4, to=2-4, loop, in=240, out=300, distance=5mm]
	\arrow["{B_1}", shift left, from=2-4, to=2-5]
	\arrow["{A_0}", shift left, bend right=10, from=2-4, to=3-3]
	\arrow["{A_1}", shift left, from=2-5, to=2-4]
	\arrow[shift left, from=2-5, to=2-6]
	\arrow[shift left, from=2-6, to=2-5]
	\arrow["{C_{r-2,0}}", from=2-6, to=2-6, loop, in=60, out=120, distance=5mm]
	\arrow["{C_{r-2,1}}", from=2-6, to=2-6, loop, in=240, out=300, distance=5mm]
	\arrow["{B_{r-2}}", shift left, bend left=30, from=2-6, to=3-7]
	\arrow["{B_0}", bend left=30, from=3-3, to=2-4]
	\arrow["{C_{00}}", from=3-3, to=3-3, loop, in=240, out=300, distance=5mm]
	\arrow["{C_{01}}", from=3-3, to=3-3, loop, in=150, out=210, distance=5mm]
	\arrow["{A_{r-1}}", shift left, from=3-3, to=3-7]
	\arrow["{A_{r-2}}", shift left,  bend right=20, from=3-7, to=2-6]
	\arrow["{B_{r-1}}", shift left, from=3-7, to=3-3]
	\arrow["{C_{r-1,1}}", from=3-7, to=3-7, loop, in=240, out=300, distance=5mm]
	\arrow["{C_{r-1,0}}", from=3-7, to=3-7, loop, in=330, out=30, distance=5mm]
\end{tikzcd}
\caption{Framed quiver of $A_{r-1}$-type.}
\label{fig-framed-quiver}
\end{figure}
where the irreducible $\bb Z_r$-representations $\rho_0,...,\rho_{r-1}$
correspond to the vertices. Let $R=\rho_0^{n_0}\oplus\cdots\oplus\rho_{r-1}^{n_{r-1}}$ and consider the following quiver representation of $\widetilde{Q}$ of dimension vector $\mathbf{v}=(n_0,...,n_{r-1},1)$:
\begin{small}
\begin{equation*}
M:=\bigoplus_{i=0}^{r-1}(\bb C_{3,4}^2\otimes\End V_i)
\oplus\left(\bb C_1^1\otimes\bigoplus_{i=0}^{r-1}\Hom(V_i,V_{i+1})\right)
\oplus\left(\bb C_2^1\otimes\bigoplus_{i=0}^{r-1}\Hom(V_i,V_{i-1})\right)
\oplus\Hom(W_0,V_0),
\end{equation*}
\end{small}
where we denote $V_i:=V_{i\mmod r}$, let $\dim V_i=n_i$
and $\dim W_0=1$. Here we let $\bb C^4:=\left\langle e_1,...,e_4\right\rangle$
and its subspaces as $\bb C_{3,4}^2=\left\langle e_3,e_4\right\rangle$ and so on.

Let $G:=\prod_{i=0}^{r-1}\GL(V_i)$ be the gauge group with the canonical action on $M$ as follows
\begin{equation}\label{action-G=GL}
    \begin{aligned}
        &(g_0,...,g_{r-1})\cdot(C_{00},C_{01},...,C_{r-1,0},C_{r-1,1},
        B_0,...,B_{r-1},A_0,...,A_{r-1},u)\\
        =&\left(g_0C_{00}g_0^{-1},g_0C_{01}g_0^{-1},...,g_{r-1}C_{r-1,0}g_{r-1}^{-1},g_{r-1}C_{r-1,1}g_{r-1}^{-1},
        g_1B_0g_0^{-1},...,\right.\\
        &\left.g_0B_{r-1}g_{r-1}^{-1},
        g_0A_0g_1^{-1},...,g_{r-1}A_{r-1}g_0^{-1},g_0u\right).
    \end{aligned}
\end{equation}
Let $U\subset M$ be an open subscheme consisting of points
such that $u$ generates the $\bb CQ$-module $(V_a,C_{bc},B_d,A_e)$
where $Q$ be the unframed quiver associated to $\widetilde{Q}$
; see Figure \ref{fig-unframed-quiver} below.
Then the $G$-action on $U$ is free; see \cite[~\S 1.2]{Sze08}
for details.

\begin{defi}
We define the \textit{non-commutative Hilbert scheme of $[\bb C^4/\bb Z_r]$} as
\[\ncHilb^{R}([\bb C^4/\bb Z_r]):=[U/G]\subset[M/G],\]
which is a smooth quasi-projective scheme.
\end{defi}
The following relation between $\Hilb^R([\bb C^4/\bb Z_r])$
and $\ncHilb^{R}([\bb C^4/\bb Z_r])$ is well-known.
\begin{lem}\label{lem-relation-nc}
We have a closed embedding $\Hilb^R([\bb C^4/\bb Z_r])\hookrightarrow\ncHilb^{R}([\bb C^4/\bb Z_r])$ defined by the following relations: for any $i\equiv 0,...,r-1\mmod r$ and $j=0,1$,
we have
\[
[C_{i0},C_{i1}]=0,\quad A_iB_i=B_{i-1}A_{i-1},\quad 
B_iC_{ij}=C_{i+1,j}B_i,\quad C_{ij}A_i=A_iC_{i+1,j}.
\]
\end{lem}

\subsubsection{Non-commutative moduli stack of sheaves \texorpdfstring{$\ncM_{R}([\bb C^4/\bb Z_r])$}{ncM R [C4/Zr]}}\label{sec-noncommu-ncMR}
Consider the following unframed quiver with the similar notation as above.
\begin{figure}[!ht]
    \begin{tikzcd}
	{Q:=} &&& {\bullet V_1} & \cdots & {\bullet V_{r-2}} \\
	&& {\bullet V_0} &&&& {\bullet V_{r-1},}
	\arrow["{C_{11}}", from=1-4, to=1-4, loop, in=60, out=120, distance=5mm]
	\arrow["{C_{10}}", from=1-4, to=1-4, loop, in=240, out=300, distance=5mm]
	\arrow["{B_1}", shift left, from=1-4, to=1-5]
	\arrow["{A_0}", shift left, bend right=10, from=1-4, to=2-3]
	\arrow["{A_1}", shift left, from=1-5, to=1-4]
	\arrow[shift left, from=1-5, to=1-6]
	\arrow[shift left, from=1-6, to=1-5]
	\arrow["{C_{r-2,0}}", from=1-6, to=1-6, loop, in=60, out=120, distance=5mm]
	\arrow["{C_{r-2,1}}", from=1-6, to=1-6, loop, in=240, out=300, distance=5mm]
	\arrow["{B_{r-2}}", shift left, bend left=30, from=1-6, to=2-7]
	\arrow["{B_0}", bend left=30, from=2-3, to=1-4]
	\arrow["{C_{00}}", from=2-3, to=2-3, loop, in=240, out=300, distance=5mm]
	\arrow["{C_{01}}", from=2-3, to=2-3, loop, in=150, out=210, distance=5mm]
	\arrow["{A_{r-1}}", shift left, from=2-3, to=2-7]
	\arrow["{A_{r-2}}", shift left, bend right=20, from=2-7, to=1-6]
	\arrow["{B_{r-1}}", shift left, from=2-7, to=2-3]
	\arrow["{C_{r-1,1}}", from=2-7, to=2-7, loop, in=240, out=300, distance=5mm]
	\arrow["{C_{r-1,0}}", from=2-7, to=2-7, loop, in=330, out=30, distance=5mm]
\end{tikzcd}
\caption{Unframed quiver with $A_{r-1}$-type.}
\label{fig-unframed-quiver}
\end{figure}

Let $R=\rho_0^{n_0}\oplus\cdots\oplus\rho_{r-1}^{n_{r-1}}$ again and consider the following quiver representation of $Q$ of dimension vector $\mathbf{v}'=(n_0,...,n_{r-1})$:
\begin{equation*}
M^{\mathrm{un}}:=\bigoplus_{i=0}^{r-1}(\bb C_{3,4}^2\otimes\End V_i)
\oplus\left(\bb C_1^1\otimes\bigoplus_{i=0}^{r-1}\Hom(V_i,V_{i+1})\right)
\oplus\left(\bb C_2^1\otimes\bigoplus_{i=0}^{r-1}\Hom(V_i,V_{i-1})\right),
\end{equation*}
where here we also denote $V_i:=V_{i\mmod r}$ and $\dim V_i=n_i$ as before.
Let the gauge group $G:=\prod_{i=0}^{r-1}\GL(V_i)$ act on $M^{\mathrm{un}}$ in a similar way to (\ref{action-G=GL}).
\begin{defi}
We define the \textit{non-commutative moduli stack of zero dimensional sheaves of $[\bb C^4/\bb Z_r]$} as
\[\ncM_{R}([\bb C^4/\bb Z_r]):=[M^{\mathrm{un}}/G]\]
which is a smooth Artin stack.
\end{defi}

We also have closed embedding $\MSh_R([\bb C^4/\bb Z_r])\hookrightarrow\ncM_R([\bb C^4/\bb Z_r])$
defined by the same relations in Lemma \ref{lem-relation-nc}.
Finally, there is a natural forgetful morphism 
\[\mathsf{nc}\mathfrak f:\ncHilb^R([\bb C^4/\bb Z_r])\to\ncM_{R}([\bb C^4/\bb Z_r])\]
by forgetting $\Hom(W_0,V_0)$.
Restricting $\mathsf{nc}\mathfrak f$ to the zero loci of the relations in Lemma \ref{lem-relation-nc},
we obtain the natural forgetful map $\mathfrak f:\Hilb^R([\bb C^4/\bb Z_r])\to\MSh_{R}([\bb C^4/\bb Z_r])$ defined in the usual way.

\subsection{Orientations of \texorpdfstring{$\Hilb^{R}([\mathbb C^4/\bb Z_r])$}{Hilb R [C4/Zr]} via non-commutative Hilbert schemes}\label{sec-ori-A}
Recall that Kool-Rennemo \cite{KR} established a method to determine the sign rules, which we summarize as follows:
\begin{itemize}
    \item Let $M$ be a quasi-projective scheme with a $3$-term symmetric obstruction theory arising from a $(-2)$-shifted symplectic derived enhancement. Suppose $M$ admits a $\bb T$-action by some algebraic torus $\bb T$ such that the obstruction theory $\bb E$ is $\bb T$-equivariant and $M^{\bb T}$ is isolated reduced. Fix a self-dual resolution
\[
\bb E\cong\{F\xrightarrow{a} E\cong E^{\vee}\xrightarrow{a^{\vee}} F^{\vee}\},
\]
where $(E,q)$ is a quadratic bundle of even rank. Given a $\bb T$-equivariant maximal isotropic subbundle $\Lambda\subset E$ inducing a canonical orientation $o$, it was proved in \cite[Prop.~2.29]{KR} that
\begin{equation}\label{eqn-virtual-class-general-points}
\begin{aligned}
  \relax [M]_{\bb T,o}^{\vir}&=\sum_{p\in M^{\bb T},\vdim\{p\}=0}(-1)^{n_{p}^{\Lambda}}\iota_{p,*}\left(\frac{1}{e_{\bb T}(F-\Lambda)|_{p}}\right),
\end{aligned}
\end{equation}
where $n_{p}^{\Lambda}\equiv \dim\left(\coker(p_{\Lambda}\circ a|_p)\right)^{\bb T\text{-}\fix}\mmod 2$ and $p_{\Lambda}$ be the projection for any splitting
$E|_p=\Lambda|_p\oplus\Lambda^{\vee}|_p$.
\end{itemize} 

In our setting, we construct a quadratic bundle over $\ncHilb^{R}([\bb C^4/\bb Z_r])$ equipped 
with a canonical maximal isotropic subbundle. This construction realizes 
$\Hilb^{R}([\bb C^4/\bb Z_r])$ as the zero locus of an isotropic section of this 
quadratic bundle. See Remark \ref{nota-2-derived-obs-ori}
for the interpretation via derived geometry.

\begin{prp}\label{prp-zero-section-Hilb-C4Zr}
Consider the following trivial bundle on $U$:
\begin{equation}
E:=\left(
    \begin{aligned}
        &\bigoplus_{k=0}^{r-1}(\Lambda^2\bb C_{34}^2\otimes\End(V_k))\oplus\bigoplus_{k=0}^{r-1}(\Lambda^2\bb C_{12}^2\otimes\End(V_k))\\
        \oplus&\bigoplus_{k=0}^{r-1}(\Lambda^2\bb C_{31}^2\otimes\Hom(V_k,V_{k+1}))\oplus\bigoplus_{k=0}^{r-1}(\Lambda^2\bb C_{32}^2\otimes\Hom(V_k,V_{k-1}))\\
        \oplus&\bigoplus_{k=0}^{r-1}(\Lambda^2\bb C_{41}^2\otimes\Hom(V_k,V_{k+1}))\oplus\bigoplus_{k=0}^{r-1}(\Lambda^2\bb C_{42}^2\otimes\Hom(V_k,V_{k-1}))
    \end{aligned}
\right)\otimes\scr{O}_U
\end{equation}
and define the quadratic form $q:E\otimes E\to\Lambda^4\bb C^4\cong\bb C$ over $E$ by
\begin{align*}
    &\sum_{k=0}^{r-1}\left(w_1^k\otimes f_1^k+w_2^k\otimes f_2^k+\cdots+w_6^k\otimes f_6^k\right)
    \otimes\sum_{k=0}^{r-1}\left(v_1^k\otimes g_1^k+v_2^k\otimes g_2^k+\cdots+v_6^k\otimes g_6^k\right)\\
    \mapsto&\sum_{k=0}^{r-1}\left(w_1^k\wedge v_2^k\cdot\tr(f_1^k\circ g_2^k)+w_2^k\wedge v_1^k\cdot\tr(f_2^k\circ g_1^k)
    +w_3^k\wedge v_6^{k+1}\cdot\tr(f_3^k\circ g_6^{k+1})\right.\\
    &\left. +w_4^k\wedge v_5^{k-1}\cdot\tr(f_4^k\circ g_5^{k-1})
    +w_5^k\wedge v_4^{k+1}\cdot\tr(f_5^k\circ g_4^{k+1})
    +w_6^k\wedge v_3^{k-1}\cdot\tr(f_6^k\circ g_3^{k-1})\right),
\end{align*}
which is nondegenerate by linear algebra.

Define a section $s\in\Gamma(U,E)$ given by
\begin{tiny}
\[\begin{tikzcd}
	{\left(e_3\otimes C_{00},e_4\otimes C_{01},...,e_3\otimes C_{r-1,0},e_4\otimes C_{r-1,1},e_1\otimes B_0,...,e_1\otimes B_{r-1},e_2\otimes A_0,...,e_2\otimes A_{r-1},u\right)} \\
	\\
	\begin{array}{c} \begin{aligned} &\sum_{k=0}^{r-1}\left((e_3\wedge e_4)\otimes(C_{k0}\circ C_{k1})+(e_4\wedge e_3)\otimes(C_{k1}\circ C_{k0})\right)+\sum_{k=0}^{r-1}\left((e_2\wedge e_1)\otimes(A_k\circ B_k)+(e_1\wedge e_2)\otimes(B_{k-1}\circ A_{k-1})\right)\\ +&\sum_{k=0}^{r-1}\left((e_1\wedge e_3)\otimes(B_k\circ C_{k0})+(e_3\wedge e_1)\otimes(C_{k+1,0}\circ B_k)\right)+\sum_{k=0}^{r-1}\left((e_3\wedge e_2)\otimes(C_{k0}\circ A_k)+(e_2\wedge e_3)\otimes(A_{k}\circ C_{k+1,0})\right)\\ +&\sum_{k=0}^{r-1}\left((e_1\wedge e_4)\otimes(B_k\circ C_{k1})+(e_4\wedge e_1)\otimes(C_{k+1,1}\circ B_k)\right)+\sum_{k=0}^{r-1}\left((e_4\wedge e_2)\otimes(C_{k1}\circ A_k)+(e_2\wedge e_4)\otimes(A_{k}\circ C_{k+1,1})\right). \end{aligned} \end{array}
	\arrow["s", maps to, from=1-1, to=3-1]
\end{tikzcd}\]
\end{tiny}
Define a subbundle $\Lambda\subset E$ given by 
\begin{equation}
\Lambda:=\left(
    \begin{aligned}
        &\bigoplus_{k=0}^{r-1}(\Lambda^2\bb C_{34}^2\otimes\End(V_k))\\
        \oplus&\bigoplus_{k=0}^{r-1}(\Lambda^2\bb C_{41}^2\otimes\Hom(V_k,V_{k+1}))\oplus\bigoplus_{k=0}^{r-1}(\Lambda^2\bb C_{42}^2\otimes\Hom(V_k,V_{k-1}))
    \end{aligned}
\right)\otimes\scr{O}_U\subset E.
\end{equation}
As all the data $(E,q,s,\Lambda)$ are $G$-equivariant,
they can descend to $\ncHilb^R([\bb C^4/\bb Z_r])$ and we will use the same notation for simplicity.
Then the section $s$ is $q$-isotropic such that $Z(s)\cong\Hilb^R([\bb C^4/\bb Z_r])$ and $\Lambda\subset E$ is a maximal isotrpic subbundle.
\end{prp}
\begin{proof}
The facts that section $s$ and $\Lambda$ are $q$-isotropic 
follow from Jacobi identity.
Moreover, since $\dim\Lambda=\frac{1}{2}\dim E$, it is maximal isotropic.
Finally, $Z(s)\cong\Hilb^R([\bb C^4/\bb Z_r])$ follows directly from Lemma \ref{lem-relation-nc}.
\end{proof}

Then the obstruction theory induced by Proposition \ref{prp-zero-section-Hilb-C4Zr} is a symmetric complex
\begin{equation}\label{eqn-obs-theory-C4Zr-2}
    T_{\ncHilb^R([\bb C^4/\bb Z_r])}\xrightarrow{ds} E\cong^qE^{\vee}\xrightarrow{ds^{\vee}} \Omega_{\ncHilb^R([\bb C^4/\bb Z_r])}.
\end{equation}
There is an orientation $o_A$ of the obstruction theory given by
choosing $\Lambda\subset E$ and $T_{\ncHilb^R([\bb C^4/\bb Z_r])}$.
By Theorem \ref{toric-vertex-thm}, we know that $\vdim\{p\}=0$ for arbitrary $p\in\Hilb^R([\bb C^4/\bb Z_r])^{\bb T}$.
Hence using Eqn.(\ref{eqn-virtual-class-general-points}), we have
\begin{equation}\label{Eqn-obs-theory-C4Zr}
\begin{aligned}
\relax    [\Hilb^{R}([\bb C^4/\bb Z_r])]_{\bb T,o_A}^{\vir}
    &=
\sum_{p\in\Hilb^R([\bb C^4/\bb Z_r])^{\bb T}}(-1)^{n_p^{\Lambda}}
\iota_{p,*}\left(\frac{1}{e_{\bb T}(T_{\ncHilb^R([\bb C^4/\bb Z_r])}-\Lambda)|_p}\right)\\
&=(-1)^{\rank\Lambda}\sum_{p\in\Hilb^R([\bb C^4/\bb Z_r])^{\bb T}}(-1)^{n_p^{\Lambda^{\vee}}}
\iota_{p,*}\left(\frac{1}{e_{\bb T}(T_{\ncHilb^R([\bb C^4/\bb Z_r])}-\Lambda^{\vee})|_p}\right),
\end{aligned}
\end{equation}
where $n_p^{\Lambda}\equiv\dim\left(\coker(p_{\Lambda}\circ ds|_p)^{\bb T\text{-}\fix}\right)\mmod2$
and $n_p^{\Lambda^{\vee}}\equiv\dim\left(\coker(p_{\Lambda^{\vee}}\circ ds|_p)^{\bb T\text{-}\fix}\right)\mmod2$ and $p_{\Lambda}$
and $p_{\Lambda^{\vee}}$ are the projections of any splitting
$E|_p=\Lambda|_p\oplus\Lambda^{\vee}|_p$. Then we have the following comparison of Eqn.~(\ref{eqn-obs-theory-C4Zr}) and Eqn.~(\ref{eqn-obs-theory-C4Zr-2}).

\begin{prp}\label{prp-vertex-T-Lambda-comparision}
For any $p_{\pi}\in\Hilb^R([\bb C^4/\bb Z_r])^{\bb T}$
correspond to the solid partition $\pi$ of color $R$, we have
\begin{equation*}
\left(T_{\ncHilb^R([\bb C^4/\bb Z_r])}-\Lambda^{\vee}\right)|_{p_{\pi}}
+\left(\Omega_{\ncHilb^R([\bb C^4/\bb Z_r])}-\Lambda\right)|_{p_{\pi}}
=v^{\DT,\bb Z_r}_{\pi}+\overline{v^{\DT,\bb Z_r}_{\pi}}.
\end{equation*}
Hence $\left(T_{\ncHilb^R([\bb C^4/\bb Z_r])}-\Lambda^{\vee}\right)|_{p_{\pi}}$ has no $\bb T$-fixed terms.
Finally, we have
\begin{equation}\label{prp-eqn-ev-eOVT}
e_{\bb T\times\bb C_m^*}\left.\left((\scr O^{[R]})^{\vee}\otimes y+\Lambda^{\vee}-T_{\ncHilb^R([\bb C^4/\bb Z_r])}\right)\right|_{p_{\pi}}
=(-1)^{\dim R+k_{\pi}}e_{\bb T\times\bb C_m^*}\left(-\widetilde{v}_{\pi}^{\DT,\bb Z_r}\right),
\end{equation}
where $k_{\pi}:=|\{(i,j,k,l)\in\pi:l\neq\min\{i,j,k\}\}|$.
\end{prp}

\begin{proof}
For any $x\in\ncHilb^R([\bb C^4/\bb Z_r])$, we have
the short exact sequence
\begin{align*}
0\to\bigoplus_{a=0}^{r-1}\End(V_a)\to&\bigoplus_{a=0}^{r-1}(\bb C_{34}^2\otimes\End(V_a))\oplus\bigoplus_{a=0}^{r-1}(\bb C_1^1\otimes\Hom(V_a,V_{a+1})\\
&\oplus\bigoplus_{a=0}^{r-1}(\bb C_2^1\otimes\Hom(V_a,V_{a-1})
\oplus\Hom(W_0,V_0)\to T_{\ncHilb^R([\bb C^4/\bb Z_r]),x}\to0.
\end{align*}
This induces
\begin{align*}
\left(T_{\ncHilb^R([\bb C^4/\bb Z_r])}-\Lambda^{\vee}\right)|_{p_{\pi}}
&=(t_3^{-1}+t_4^{-1}-t_1^{-1}t_2^{-1}-1)\sum_{a=0}^{r-1}Z_{\pi}^{(a)}\overline{Z_{\pi}^{(a)}}\\
&+
(t_1^{-1}-t_1^{-1}t_3^{-1})\sum_{a-b\equiv 1\mmod r}Z_{\pi}^{(a)}\overline{Z_{\pi}^{(b)}}\\
&+(t_2^{-1}-t_2^{-1}t_3^{-1})\sum_{a-b\equiv -1\mmod r}Z_{\pi}^{(a)}\overline{Z_{\pi}^{(b)}}+Z_{\pi}^{(0)}
\end{align*}
as K-theoretic classes.
By direct calculation, we have
\[
\left(T_{\ncHilb^R([\bb C^4/\bb Z_r])}-\Lambda^{\vee}\right)|_{p_{\pi}}
-v^{\DT,\bb Z_r}_{\pi}=(t_4^{-1}-t_4)\sum_{a=0}^{r-1}Z_{\pi}^{(a)}\overline{Z_{\pi}^{(a)}},
\]
which induce 
\[\left(T_{\ncHilb^R([\bb C^4/\bb Z_r])}-\Lambda^{\vee}\right)|_{p_{\pi}}
+\left(\Omega_{\ncHilb^R([\bb C^4/\bb Z_r])}-\Lambda\right)|_{p_{\pi}}
=v^{\DT,\bb Z_r}_{\pi}+\overline{v^{\DT,\bb Z_r}_{\pi}}.\]

By Theorem \ref{toric-vertex-thm}, $v_{\pi}^{\DT}$ has no $\bb T$-fixed terms. Then 
\[v^{\DT,\bb Z_r}=(v^{\DT})^{\bb Z_r\text{-}\fix}
\text{ and }\left(T_{\ncHilb^R([\bb C^4/\bb Z_r])}-\Lambda^{\vee}\right)|_{p_{\pi}}\]
also have no $\bb T$-fixed terms; see also \cite[Lem.~5.9]{CKM23}.
Hence, we have 
\[e_{\bb T\times\bb C_m^*}\left.\left((\scr O^{[R]})^{\vee}\otimes y+\Lambda^{\vee}-T_{\ncHilb^R([\bb C^4/\bb Z_r])}\right)\right|_{p_{\pi}}
=(-1)^{\lambda}e_{\bb T\times\bb C_m^*}\left(-\widetilde{v}_{\pi}^{\DT,\bb Z_r}\right),\]
where $\lambda=\dim\left(t_4^{-1}\sum_{a=0}^{r-1}Z_{\pi}^{(a)}\overline{Z_{\pi}^{(a)}}\right)^{\bb T\text{-}\mov}$.
This is given by
\begin{equation}\label{eqn-k-pi-calculation}
\begin{aligned}
\lambda&=\sum_{a=0}^{r-1}n_a^2-\dim\left(t_4^{-1}\sum_{a=0}^{r-1}Z_{\pi}^{(a)}\overline{Z_{\pi}^{(a)}}\right)^{\bb T\text{-}\fix}\\
    &\equiv\dim R-\dim\left(t_4^{-1}Z_{\pi}\overline{Z_{\pi}}\right)^{(\bb T\times\bb Z_r)\text{-}\fix}
    =\dim R-\dim\left(t_4^{-1}Z_{\pi}\overline{Z_{\pi}}\right)^{\bb T\text{-}\fix}\\
    &\equiv\dim R-k_{\pi}\equiv\dim R+k_{\pi}\mmod 2,
\end{aligned}
\end{equation}
where the third equality follows from the fact that 
\[\left(t_4^{-1}Z_{\pi}\overline{Z_{\pi}}\right)^{(\bb T\times\bb Z_r)\text{-}\fix}=\left(t_4^{-1}Z_{\pi}\overline{Z_{\pi}}\right)^{\bb T\text{-}\fix}\] 
and the fourth equality follows from \cite[Lem.~4.3]{KR}. 
\end{proof}

\begin{nota}\label{nota-2-derived-obs-ori}
Note that this obstruction theory is given by another derived enhancement 
$\mathsf{Z}^{\mathrm{symp}}_{\ncHilb^R([\bb C^4/\bb Z_r])}(E,s)$ of $Z(s)\cong\Hilb^R([\bb C^4/\bb Z_r])$, which we call it the 
derived symplectic zero locus of isotropic section $s$ (see \cite[\S~3.2,~Cor.~4.2.2]{park24} for details). 
Here we note that the derived Hilbert scheme 
$\dHilb^R([\bb C^4/\bb Z_r])$ is 
believed (but still conjectured) to be isomorphic to $\mathsf{Z}^{\mathrm{symp}}_{\ncHilb^R([\bb C^4/\bb Z_r])}(E,s)$.
\end{nota}

Therefore, we can deduce the following formula.
\begin{prp}
We have 
\[
\int_{[\Hilb^{R}([\bb C^4/\bb Z_r])]^{\vir}_{\bb T,o_A}}e_{\bb T\times\bb C^*_m}((\scr O^{[R]})^{\vee}\otimes e^m)
    =\sum_{|\pi|_{\bb Z_r}=R}(-1)^{n_{p_{\pi}}^{\Lambda^{\vee}}+k_{\pi}}
    e_{\bb T\times\bb C_m^*}\left(-\widetilde{v}_{\pi}^{\DT,\bb Z_r}\right).
\]
\end{prp}
\begin{proof}
By Eqn. (\ref{Eqn-obs-theory-C4Zr}), we have 
\begin{align*}
 [\Hilb^{R}([\bb C^4/\bb Z_r])]_{\bb T,o_A}^{\vir}
&=(-1)^{\dim R}\sum_{p\in\Hilb^R([\bb C^4/\bb Z_r])^{\bb T}}(-1)^{n_p^{\Lambda^{\vee}}}
\iota_{p,*}\left(\frac{1}{e_{\bb T}(T_{\ncHilb^R([\bb C^4/\bb Z_r])}-\Lambda^{\vee})|_p}\right),
\end{align*}
because
$\rank\Lambda=\sum_{a=0}^{r-1}n_a^2+2(n_0n_1+\cdots+n_{r-1}n_0)
\equiv\sum_{a=0}^{r-1}n_a^2\equiv\dim R\mmod2$.
Hence
\begin{align*}
&\int_{[\Hilb^{R}([\bb C^4/\bb Z_r])]^{\vir}_{\bb T,o_A}}e_{\bb T\times\bb C^*_m}((\scr O^{[R]})^{\vee}\otimes e^m)\\
    =&(-1)^{\dim R}\sum_{|\pi|_{\bb Z_r}=R}(-1)^{n_{p_{\pi}}^{\Lambda^{\vee}}}
\frac{e_{\bb T\times\bb C^*_m}((\scr O^{[R]})^{\vee}\otimes y)|_{p_{\pi}}}{e_{\bb T}\left.\left(T_{\ncHilb^R([\bb C^4/\bb Z_r])}-\Lambda^{\vee}\right)\right|_{p_{\pi}
}}\\
    =&(-1)^{\dim R}\sum_{|\pi|_{\bb Z_r}=R}(-1)^{n_{p_{\pi}}^{\Lambda^{\vee}}}
    \left.e_{\bb T\times\bb C^*_m}\left((\scr O^{[R]})^{\vee}\otimes y+\Lambda^{\vee}-T_{\ncHilb^R([\bb C^4/\bb Z_r])}\right)\right|_{p_{\pi}}\label{eqn-ch-to-barket}
    \\
    =&\sum_{|\pi|_{\bb Z_r}=R}
    (-1)^{n_{p_{\pi}}^{\Lambda^{\vee}}+k_{\pi}}e_{\bb T\times\bb C^*_m}\left(-\widetilde{v}_{\pi}^{\DT,\bb Z_r}\right),
\end{align*}
where the final equality follows from Eqn. (\ref{prp-eqn-ev-eOVT}) in Proposition \ref{prp-vertex-T-Lambda-comparision}.
\end{proof}

Finally, we need to compute $n_p^{\Lambda^{\vee}}+k_{\pi}\mmod2$. Firstly, we need the following lemma.

\begin{lem}[{\cite[Lem.~4.6]{KR}}]\label{lem-KR-4.9}
Denote by $p:\bb C^4\to\bb C^3$ the projection onto the first three coordinates. For any $\bb T$-fixed $0$-dimensional subscheme 
$Z\subset\bb C^4$ corresponding to a
solid partition $\pi$ we have\footnote{Note that $\bb T\subset(\bb C^*)^4$ defined by relation $t_1t_2t_3t_4=1$. Projecting $\bb T$ onto the first three coordinates, we get $\bb T\cong(\bb C^*)^3$ which acts in the standard way on $\bb C^3$.}
\begin{enumerate}
    \item $\dim\Hom(p_*\scr O_Z,p_*\scr O_Z)^{\bb T\text{-}\fix}=|(i,i,i,j)\in\pi:j\geq i|$;
    \item $\dim\Ext^1(p_*\scr O_Z,p_*\scr O_Z)^{\bb T\text{-}\fix}=\dim\Ext^2(p_*\scr O_Z,p_*\scr O_Z)^{\bb T\text{-}\fix}=0$;
    \item $\dim\Ext^3(p_*\scr O_Z,p_*\scr O_Z)^{\bb T\text{-}\fix}=|(i,i,i,j)\in\pi:j> i|$.
\end{enumerate}
\end{lem}

\begin{prp}
In this case, we have
\[n_{p_{\pi}}^{\Lambda^{\vee}}+k_{\pi}\equiv\mu_{\pi}\mmod2,\]
where $\mu_{\pi}=|\{(i,i,i,j)\in\pi:j>i\}|$.
\end{prp}
\begin{proof}
This is a generalization of \cite[Prop.~4.5]{KR}.
Recall that
\[n_{p_{\pi}}^{\Lambda^{\vee}}\equiv\dim\left(\coker(p_{\Lambda^{\vee}}\circ ds|_{p_{\pi}})^{\bb T\text{-}\fix}\right)\mmod2,\]
where $0\to\Lambda\to E\xrightarrow{p_{\Lambda^{\vee}}}\Lambda^{\vee}\to0$
and $ds:T_{\ncHilb^R([\bb C^4/\bb Z_r])}\to E$.
Let 
\begin{small}
    \[
\left(e_3\otimes C_{00}^{\pi},e_4\otimes C_{01}^{\pi},...,e_3\otimes C_{r-1,0}^{\pi},e_4\otimes C_{r-1,1}^{\pi},e_1\otimes B_0^{\pi},...,e_1\otimes B_{r-1}^{\pi},e_2\otimes A_0^{\pi},...,e_2\otimes A_{r-1}^{\pi},u^{\pi}\right)
\]
\end{small}
be the representation of $p_{\pi}$.
Consider the composition 
\begin{align*}
&\bigoplus_{a=0}^{r-1}(\bb C_{34}^2\otimes\End(V_a))\oplus\bigoplus_{a=0}^{r-1}(\bb C_1^1\otimes\Hom(V_a,V_{a+1})\\
&\oplus\bigoplus_{a=0}^{r-1}(\bb C_2^1\otimes\Hom(V_a,V_{a-1})
\oplus\Hom(W_0,V_0)\twoheadrightarrow T_{\ncHilb^R([\bb C^4/\bb Z_r]),p_{\pi}}\\
&\xrightarrow{ds|_{p_{\pi}}}E|_{p_{\pi}}\xrightarrow{}\Lambda^{\vee}|_{p_{\pi}},
\end{align*}
where the first morphism is surjective. In order to consider the cokernel of $p_{\Lambda^{\vee}}\circ ds|_{p_{\pi}}$,
we just need to consider the cokernel $\coker\Theta$ of the induced morphism $\Theta$:
\begin{small}
\[\begin{tikzcd}
	{\bigoplus_{i=0}^{r-1}(\bb C_{3}^1\otimes\End V_i) \oplus\left(\bb C_1^1\otimes\bigoplus_{i=0}^{r-1}\Hom(V_i,V_{i+1})\right) \oplus\left(\bb C_2^1\otimes\bigoplus_{i=0}^{r-1}\Hom(V_i,V_{i-1})\right)} \\
	\\
	{\bigoplus_{k=0}^{r-1}(\Lambda^2\bb C_{12}^2\otimes\End(V_k))         \oplus\bigoplus_{k=0}^{r-1}(\Lambda^2\bb C_{31}^2\otimes\Hom(V_k,V_{k+1}))\oplus\bigoplus_{k=0}^{r-1}(\Lambda^2\bb C_{32}^2\otimes\Hom(V_k,V_{k-1}))}
	\arrow["\Theta", from=1-1, to=3-1]
\end{tikzcd}\]
\end{small}
given by 
\begin{small}
\[\begin{tikzcd}
	{(e_3\otimes C_{00},...,e_3\otimes C_{r-1,0},e_1\otimes B_0,...,e_1\otimes B_{r-1},e_2\otimes A_0,...,e_2\otimes A_{r-1})} \\
	\\
	\begin{array}{c} \begin{aligned}     &\sum_{k=0}^{r-1}\left((e_2\wedge e_1)\otimes(A_k^{\pi}\circ B_k-B_{k-1}\circ A_{k-1}^{\pi})+(e_1\wedge e_2)\otimes(B_{k-1}^{\pi}\circ A_{k-1}-A_k\circ B_k^{\pi})\right)\\      +&\sum_{k=0}^{r-1}\left((e_1\wedge e_3)\otimes(B_k^{\pi}\circ C_{k0}-C_{k+1,0}\circ B_k^{\pi})+(e_3\wedge e_1)\otimes(C_{k+1,0}^{\pi}\circ B_k-B_k\circ C_{k0}^{\pi})\right)\\     +&\sum_{k=0}^{r-1}\left((e_3\wedge e_2)\otimes(C_{k0}^{\pi}\circ A_k-A_{k}\circ C_{k+1,0}^{\pi})+(e_2\wedge e_3)\otimes(A_{k}^{\pi}\circ C_{k+1,0}-C_{k0}\circ A_k^{\pi})\right). \end{aligned} \end{array}
	\arrow["\Theta", maps to, from=1-1, to=3-1]
\end{tikzcd}\]
\end{small}

Consider the following cartesian diagram
\[\begin{tikzcd}
	{\bb C^4} && {\bb C^3} \\
	{[\bb C^4/\bb Z_r]} && {[\bb C^3/\bb Z_r]}
	\arrow["p", from=1-1, to=1-3]
	\arrow[two heads, from=1-1, to=2-1]
	\arrow["\square"{description}, draw=none, from=1-1, to=2-3]
	\arrow[two heads, from=1-3, to=2-3]
	\arrow["P", from=2-1, to=2-3]
\end{tikzcd}\]
by projecting to the first three coordinates,
then 
\[\RR\Hom_{[\bb C^3/\bb Z_r]}(P_*\scr O_{p_{\pi}},P_*\scr O_{p_{\pi}})
=\RR\Hom_{\bb C^3}(p_*\scr O_{p_{\pi}},p_*\scr O_{p_{\pi}})^{\bb Z_r\text{-}\fix}.\]
Write $V:=\oplus_{a=0}^{r-1}\rho_aV_a$, then $\RR\Hom_{\bb C^3}(p_*\scr O_{p_{\pi}},p_*\scr O_{p_{\pi}})$
can be computed by the following complex
\[\End(V)\to\bb C^3\otimes\End(V)\to\Lambda^2\bb C^3\otimes\End(V)
\to\Lambda^3\bb C^3\otimes\End(V).\]
As $\bb Z_r$ is a linear reductive group, then $\RR\Hom_{[\bb C^3/\bb Z_r]}(P_*\scr O_{p_{\pi}},P_*\scr O_{p_{\pi}})$ can be computed by 
the following complex
\[\End(V)^{\bb Z_r\text{-}\fix}\to(\bb C^3\otimes\End(V))^{\bb Z_r\text{-}\fix}\xrightarrow{\theta_1}(\Lambda^2\bb C^3\otimes\End(V))^{\bb Z_r\text{-}\fix}
\xrightarrow{\theta_2}(\Lambda^3\bb C^3\otimes\End(V))^{\bb Z_r\text{-}\fix}.\]
Observe that $\theta_1=\Theta$.
Now we have
\[\Ext^2_{[\bb C^3/\bb Z_r]}(P_*\scr O_{p_{\pi}},P_*\scr O_{p_{\pi}})
^{\bb T\text{-}\fix}=\Ext^2_{\bb C^3}(p_*\scr O_{p_{\pi}},p_*\scr O_{p_{\pi}})
^{(\bb T\times\bb Z_r)\text{-}\fix}=0\]
by Lemma \ref{lem-KR-4.9}(2).
Hence $\coker(\theta_1)^{\bb T\text{-}\fix}=\operatorname{Im}(\theta_2)^{\bb T\text{-}\fix}$. Now, we have
\[\dim\operatorname{Im}(\theta_2)^{\bb T\text{-}\fix}
=\dim(\Lambda^3\bb C^3\otimes\End(V))^{(\bb T\times\bb Z_r)\text{-}\fix}-\dim\Ext^3_{\bb C^3}(p_*\scr O_{p_{\pi}},p_*\scr O_{p_{\pi}})
^{(\bb T\times\bb Z_r)\text{-}\fix}.\]
By the same reason in Eqn. (\ref{eqn-k-pi-calculation}), we have 
\begin{align*}
    &\dim(\Lambda^3\bb C^3\otimes\End(V))^{(\bb T\times\bb Z_r)\text{-}\fix}=\dim\left(t_4\sum_{a=0}^{r-1}Z_{\pi}^{(a)}\overline{Z_{\pi}^{(a)}}\right)^{\bb T\text{-}\fix}\\
    =&\dim\left(t_4^{-1}\sum_{a=0}^{r-1}Z_{\pi}^{(a)}\overline{Z_{\pi}^{(a)}}\right)^{\bb T\text{-}\fix}
    \equiv k_{\pi}\mmod2.
\end{align*}
Similarly, we have
\begin{align*}
    \dim\Ext^3_{\bb C^3}(p_*\scr O_{p_{\pi}},p_*\scr O_{p_{\pi}})
^{(\bb T\times\bb Z_r)\text{-}\fix}
=\dim\Ext^3_{\bb C^3}(p_*\scr O_{p_{\pi}},p_*\scr O_{p_{\pi}})
^{\bb T\text{-}\fix}=\mu_{\pi}
\end{align*}
by Lemma \ref{lem-KR-4.9}(3) again.
This shows
$\dim\operatorname{Im}(\theta_2)^{\bb T\text{-}\fix}
\equiv k_{\pi}-\mu_{\pi}\mmod2$ and hence
\[n_{p_{\pi}}^{\Lambda^{\vee}}+k_{\pi}\equiv\mu_{\pi}\mmod2.\]
This gives the result.
\end{proof}

In conclusion, we prove the following theorem.

\begin{thm}\label{thm-ori-A-sign-rule}
For orientation $o_A$, we have
\[
\int_{[\Hilb^{R}([\bb C^4/\bb Z_r])]^{\vir}_{\bb T,o_A}}e_{\bb T\times\bb C^*_m}((\scr O^{[R]})^{\vee}\otimes e^m)
    =\sum_{|\pi|_{\bb Z_r}=R}(-1)^{\mu_{\pi}}
    e_{\bb T\times\bb C_m^*}\left(-\widetilde{v}_{\pi}^{\DT,\bb Z_r}\right),
\]
where $\mu_{\pi}=|\{(i,i,i,j)\in\pi:j>i\}|$.
\end{thm}

\subsection{Comparison of two orientations of \texorpdfstring{$\Hilb^{R}([\mathbb C^4/\bb Z_r])$}{Hilb R [C4/Zr]}}\label{section-compare-oris}
Actually, we have $[\bb C^4/\bb Z_r]=[\bb C^3/\bb Z_r]\times\bb C
\cong\Tot(\omega_{[\bb C^3/\bb Z_r]})$ by choosing the first three coordinates. Hence, by Lemma \ref{lem-ori-identification} and in particular, Eqn. (\ref{eqn-isotropic-quotient-orientation}), this gives another orientation which we call it $o_B$.
The aim here is to compare $o_A$ and $o_B$ and to give the sign rules of $o_B$,
similar to \cite[Appendix~B]{CZZ24}. Here we fix $V=\bigoplus_{a=0}^{r-1}\rho_aV_a$.

\begin{lem}\label{lem-ori-A-complex}
Consider the following diagram of vector bundles on $M^{\mathrm{un}}$:
\begin{tiny}
\[\begin{tikzcd}
	{\End(V)^{\bb Z_r\text{-}\fix}} & {(\bb C^4\otimes\End(V))^{\bb Z_r\text{-}\fix}} & {(\Lambda^2\bb C^4\otimes\End(V))^{\bb Z_r\text{-}\fix}} & {(\Lambda^3\bb C^4\otimes\End(V))^{\bb Z_r\text{-}\fix}} & {(\Lambda^4\bb C^4\otimes\End(V))^{\bb Z_r\text{-}\fix}} \\
	{\End(V)^{\bb Z_r\text{-}\fix}} & {(\bb C^4\otimes\End(V))^{\bb Z_r\text{-}\fix}} & {(\Lambda^2\bb C^3\otimes\End(V))^{\bb Z_r\text{-}\fix}} & 0 & 0,
	\arrow[from=1-1, to=1-2]
	\arrow["{=}"', from=1-1, to=2-1]
	\arrow[from=1-2, to=1-3]
	\arrow["{=}"', from=1-2, to=2-2]
	\arrow[from=1-3, to=1-4]
	\arrow[from=1-3, to=2-3]
	\arrow[from=1-4, to=1-5]
	\arrow[from=1-4, to=2-4]
	\arrow[from=1-5, to=2-5]
	\arrow[from=2-1, to=2-2]
	\arrow[from=2-2, to=2-3]
	\arrow[from=2-3, to=2-4]
	\arrow[from=2-4, to=2-5]
\end{tikzcd}\]
\end{tiny}
which is $G$-equivariant and descends to a diagram of vector bundles on $\ncM_R([\bb C^4/\bb Z_r])$ defined in $\S$\ref{sec-noncommu-ncMR}.
Its restriction to $\MSh_R([\bb C^4/\bb Z_r])$ is an isotropic quotient of complexes and induces an orientation.
The pullback of this orientation to $\Hilb^R([\bb C^4/\bb Z_r])$ via $\mathfrak f$ is given by $o_C=(-1)^{\dim R}o_A$.
\end{lem}
\begin{proof}
By the construction of $\S$\ref{sec-ori-A}, as in Eqn. (\ref{Eqn-obs-theory-C4Zr})
and the first exact sequence in Proposition \ref{prp-vertex-T-Lambda-comparision}, the differences between $o_A$ and $o_C$ are
\begin{small}
\[
\End(V)^{\bb Z_r\text{-}\fix}=\left((\Lambda^4\bb C^4\otimes\End(V))^{\bb Z_r\text{-}\fix}\right)^*, (\bb C^4\otimes\End(V))^{\bb Z_r\text{-}\fix}=\left((\Lambda^3\bb C^4\otimes\End(V))^{\bb Z_r\text{-}\fix}\right)^*.
\] 
\end{small}
Therefore, by the statement in the first paragraph of \cite[Prop.~B.2]{CZZ24},
this orientation is $(-1)^{5\dim R}o_A=(-1)^{\dim R}o_A$. Note that this is the isotropic quotient of the tangent complex rather than the cotangent complex.
\end{proof}

\begin{lem}\label{lem-ori-B-complex}
Consider the following diagram of vector bundles on $M^{\mathrm{un}}$:
\begin{equation}
\begin{tiny}
\begin{tikzcd}
	{\End(V)^{\bb Z_r\text{-}\fix}} & {(\bb C^4\otimes\End(V))^{\bb Z_r\text{-}\fix}} & {(\Lambda^2\bb C^4\otimes\End(V))^{\bb Z_r\text{-}\fix}} & {(\Lambda^3\bb C^4\otimes\End(V))^{\bb Z_r\text{-}\fix}} & {(\Lambda^4\bb C^4\otimes\End(V))^{\bb Z_r\text{-}\fix}} \\
	{\End(V)^{\bb Z_r\text{-}\fix}} & {(\bb C^3\otimes\End(V))^{\bb Z_r\text{-}\fix}} & {(\Lambda^2\bb C^3\otimes\End(V))^{\bb Z_r\text{-}\fix}} & {(\Lambda^3\bb C^3\otimes\End(V))^{\bb Z_r\text{-}\fix}} & 0,
	\arrow[from=1-1, to=1-2]
	\arrow["{=}"', from=1-1, to=2-1]
	\arrow[from=1-2, to=1-3]
	\arrow[ from=1-2, to=2-2]
	\arrow[from=1-3, to=1-4]
	\arrow[from=1-3, to=2-3]
	\arrow[from=1-4, to=1-5]
	\arrow[from=1-4, to=2-4]
	\arrow[from=1-5, to=2-5]
	\arrow[from=2-1, to=2-2]
	\arrow[from=2-2, to=2-3]
	\arrow[from=2-3, to=2-4]
	\arrow[from=2-4, to=2-5]
\end{tikzcd}
\end{tiny}
\end{equation}
which is $G$-equivariant and descends to a diagram of vector bundles on $\ncM_R([\bb C^4/\bb Z_r])$.
Its restriction to $\MSh_R([\bb C^4/\bb Z_r])$ is an isotropic quotient of complexes and induces an orientation.
The pullback of this orientation to $\Hilb^R([\bb C^4/\bb Z_r])$ via $\mathfrak f$ is given by $o_B$.
\end{lem}
\begin{proof}
This follows directly from the $\bb Z_r$-fixed version of \cite[Lem.~B.1]{CZZ24}.
\end{proof}

\begin{coro}\label{coro-sign-rule-oB}
For orientation $o_B$, we have $o_A=o_B=(-1)^{\dim R}o_C=:o$. Hence
\[
\int_{[\Hilb^{R}([\bb C^4/\bb Z_r])]^{\vir}_{\bb T,o_B=o}}e_{\bb T\times\bb C^*_m}((\scr O^{[R]})^{\vee}\otimes e^m)
    =\sum_{|\pi|_{\bb Z_r}=R}(-1)^{\mu_{\pi}}
    e_{\bb T\times\bb C_m^*}\left(-\widetilde{v}_{\pi}^{\DT,\bb Z_r}\right),
\]
where $\mu_{\pi}=|\{(i,i,i,j)\in\pi:j>i\}|$.
\end{coro}
\begin{proof}
By the resolutions of Lemma \ref{lem-ori-A-complex} and Lemma \ref{lem-ori-B-complex}, the difference between $o_B$ and $o_C$
is 
\[(\bb C_4^1\otimes\End(V))^{\bb Z_r\text{-}\fix}=\left((\Lambda^3\bb C^3\otimes\End(V))^{\bb Z_r\text{-}\fix}\right)^*.\]
By the statement in the first paragraph of \cite[Prop.~B.2]{CZZ24}
and Lemma \ref{lem-ori-A-complex} again, the difference between sign rules of $o_A$ and $o_B$ is given by
\[
\dim R+\dim(\bb C_4^1\otimes\End(V))^{\bb Z_r\text{-}\fix}
\equiv2(n_0+\cdots+n_{r-1})\equiv0\mmod 2.
\]
Therefore, the result follows from Theorem \ref{thm-ori-A-sign-rule}.
\end{proof}

\begin{nota}
The constructions and the method of comparing orientations in this section can be generalized to $[\bb C^4/G]$
for any finite subgroup $G\subset\bb T$. The comparison sign depends on $G$ and $R$.
\end{nota}

\section{\texorpdfstring{$\mathsf{DT}_4$}{DT4}-invariants of \texorpdfstring{$[\mathbb C^4/\bb Z_r]$}{[C4/Zr]}}\label{sect on DT of C4/Zr}
In this section, we will compute $\DTC_{[\bb C^4/\bb Z_r],\scr L,o}$ using the results of previous sections and the orbifold degeneration formula to deduce several pole analysis, which is the four-fold analogue of \cite[Thm.~5.3]{Zhou18-2}. Here we assume $r\geq2$ since the case of $\bb C^4$ has been proved.

\subsection{Relative and rubber \texorpdfstring{$\mathsf{DT}_4$}{DT4} virtual classes}
One of the main techniques for computing the invariant
$\DTC_{[\bb C^4/\bb Z_r],\scr L,o}$ is the expanded and rubber geometry.
The foundation of this method was first developed by Li-Wu \cite{Li01,LW15}.
In the literature, these constructions are used to establish a degeneration formula for Donaldson-Thomas invariants with several applications; see \cite{MNOP2,OP10,Zhou18,CZZ24} for details.

We note that the expanded geometry for the Deligne-Mumford stacks has been established in \cite{Zhou18}. The rubber geometry of the log CY pair $(\Delta,D_{0}\sqcup D_{\infty})$ is derived from a construction similar to \cite[~\S 4]{Zhou18-2}, where
$\Delta=[\bb C^3/\bb Z_r]\times\mathbb P^1$, $D_{0}=[\bb C^3/\bb Z_r]\times\{0\}$ and $D_{\infty}=[\bb C^3/\bb Z_r]\times\{\infty\}$.
Here we will use the techniques of expanded and rubber geometry
follows from the similar notation as in \cite{Zhou18-2,CZZ24}. 

\begin{defi}\label{defi-spaces-T-actions}
Here we state some spaces with torus actions.
\begin{enumerate}
    \item (Relative $\mathsf{DT}_4$ invariants) We define $X:=[\bb C^2/\bb Z_r]\times\Tot(\scr O_{\bb P^{1}}(-1))$ and $D_{\infty}=[\bb C^2/\bb Z_r]\times\{\infty\}\times\bb C$ with $(\bb C^*)^4$-action on $X$ with the following tangent weights at fixed points
    \[\text{at }0:-s_1,-s_2,-s_3,-s_4,\quad\text{at }\infty:-s_1,-s_2,s_3,-(s_3+s_4).\]
    We take a $3$-dim torus $\bb T\subset (\bb C^*)^4$ similar to \cite[~\S 6.3~(1)]{CZZ24} with equivariant parameters of relation $s_1+s_2+s_3+s_4=0$.
    The orientation (see \cite[Def.~4.5,~Thm.~4.6]{CZZ24}) of shifted Lagrangian is constructed from
\[X=\Tot(\omega_{[\bb C^2/\bb Z_r]\times\bb P^1}([\bb C^2/\bb Z_r]\times\{\infty\})),\quad
D_{\infty}=\Tot(\omega_{[\bb C^2/\bb Z_r]\times\{\infty\}}).\]
In this case, we have the stack of expanded pairs $\mcal A^R$ with the universal object $\mcal R^R$ and the Hilbert stack $\Hilb^R(X,D_{\infty})$.

    \item (Rubber invariants) We define $\Delta:=[\bb C^2/\bb Z_r]\times\mathbb P^1\times\bb C$ and $D_0=[\bb C^2/\bb Z_r]\times\{0\}\times\bb C$
    and $D_{\infty}=[\bb C^2/\bb Z_r]\times\{\infty\}\times\bb C$
    where the $\bb C^*$-action on $\bb P^1$ is considered as automorphism in the relative Hilbert stack. 
    In this case, $\bb T\cong(\bb C^*)^2$ is the CY torus on the fiber $[\bb C^2/\bb Z_r]\times\bb C$, whose equivariant parameters satisfy $s_1+s_2+s_4=0$. Finally, the orientation is constructed from
\[\Delta=\Tot(\omega_{[\bb C^2/\bb Z_r]\times\bb P^1}([\bb C^2/\bb Z_r]\times\{0,\infty\})),\quad
D_0\sqcup D_{\infty}=\Tot(\omega_{[\bb C^2/\bb Z_r]\times\{0,\infty\}}).\]
In this case, we have $\mcal A^{\sim}\cong\mathfrak M^{\text{ss}}_{0,2}$
and $\mcal R^{\sim}\cong\mcal C_{0,2}\times D$ and Hilbert stacks
$\Hilb^{\sim,R}(\Delta,D_0\sqcup D_{\infty})$
where $\mathfrak M^{\text{ss}}_{0,2}$ is the moduli stack of semistable curves of genus zero with two marked points
and $\mcal C_{0,2}$ is its universal curve; see \cite[~\S A.1--A.2]{CZZ24}.
\end{enumerate}
\end{defi}

\begin{nota}
By \cite[~\S 3.2]{Zhou18-2}, we have
\[F_0K_{\text{cpt}}^{\text{num}}(\Delta)=\left\langle [\rho_i\otimes\scr{O}_{p}]:
p\in\mathbf B\bb Z_r\times\{0\}\times\bb P^1,0\leq i\leq r-1
\right\rangle\cong\text{Rep}(\bb Z_r).\]
This holds similarly for $X$.
\end{nota}

The restriction maps of expanded and rubber pairs
\[r_{\mcal D\to\mcal R}:\Hilb^{R}(X, D_{\infty})\to\mcal A^{R}\]
and
\[r_{\mcal D\to\mcal R^{\sim}}:\Hilb^{\sim,R}(\Delta,D_0\sqcup D_{\infty})\to\mcal A^{\sim,R}\]
admit canonical $\bb T$-equivariant isotropic symmetric obstruction theories. These constructions induce $\bb T$-equivariant virtual pullbacks $\sqrt{r_{\mcal D\to\mcal R}^!}$
and $\sqrt{r_{\mcal D\to\mcal R^{\sim}}^!}$.

\begin{defi}\label{defi-relativeDT4-rubber-virtual-classes}
We define the relative $\mathsf{DT}_4$ and rubber virtual classes as follows.
\begin{enumerate}
    \item We define \textit{relative $\mathsf{DT}_4$ virtual class} as
    \[
    [\Hilb^R(X,D_{\infty})]^{\vir}_{\bb T}:=\sqrt{r_{\mcal D\to\mcal R}^!}[\mcal A^R].
    \]
    \item Let $\pi_{\mcal Z^{\sim,R}}:\mcal Z^{\sim,R}\to\Hilb^{\sim,R}(\Delta,D_0\sqcup D_{\infty})$
be the proper flat universal projection. We define the \textit{rubber virtual classes} by
\begin{align}
    &[\Hilb^{\sim,R}(\Delta,D_0\sqcup D_{\infty})]_{\bb T}^{\vir}:=\sqrt{r_{\mcal D\to\mcal R^{\sim}}^!}[\mcal A^{\sim,R}],\\
    &[\mcal Z^{\sim,R}]^{\vir}:=\pi_{\mcal Z^{\sim,R}}^*[\Hilb^{\sim,R}(\Delta,D_0\sqcup D_{\infty})]_{\bb T}^{\vir}.
\end{align}
\end{enumerate}
\end{defi}

\begin{nota}
Here the isotropic condition follows 
from the fact that these admit exact $(-2)$-shifted symplectic derived enhancements as defined in 
\cite{CZZ24} since we only consider Hilbert scheme of points. 

However, even for the case of curve 
counting of orbifolds,
the isotropy conditions of symmetric obstruction theories still hold
similar as \cite[Prop.~5.6,~5.11,~Appendix~C]{CZZ24}, which using the theory of \cite{park24}, via the following modifications:
\begin{itemize}
    \item Note that \cite[Eqn.~(C.4)]{CZZ24} holds since $M_b$ is a scheme; 
    see also \cite{OS03}. Note also that $H^4(\mathcal X_b,\mathcal X_b\backslash Z)=0$ in \cite{CZZ24} 
    holds for the smooth 
    pairs of schemes. Hence, one can just take the invariant part to deduce the 
    case of quotient stacks. See the proof of the second paragraph of 
    \cite[Prop.~36]{Behrend02}. Then we need to consider the comparison of 
    ``$\text{HP}_Z^k(\mcal X_b,\mcal X_b\backslash Z)\cong H_{\text{sing}}^{k+2p}(\mcal X_b,\mcal X_b\backslash Z)$''.
  \item By \cite[Cor.~3.4]{BMS19}, we have the faithfully flat descent of the period cyclic homology. Hence, by this and \cite[Thm.~2.2]{Emm95} we have $\text{HP}=H_{dR}^*$ as in \cite[Eqn.~(C.3)]{CZZ24}. 
  \item Finally, we need to consider the comparison of $H_{dR}^*$
and the singular cohomology of the underlying orbifolds. We need to follow the proof in \cite[Thm.~IV.2.1]{Har75}. Our Deligne-Mumford stack has the resolution of singularities which follows from the equivariant resolution of singularity by the finite group action as in \cite{AW97}. Then the whole proof works for our case. See also \cite{Behrend02}.
\end{itemize}
\end{nota}

\subsection{Degeneration formula and localization formula}\label{sec_degeneration_formula}
The degeneration formula given by the restriction of the universal family $\mcal Z^{\sim,R}$ to some divisor of the rubber family $\mcal R^{\sim,R}$ to make the virtual classes decomposed. Here we will recall some descriptions of the rubber family first.

A map $S \rightarrow \mathcal{R}^\sim$ is equivalent to a pair $(\mathcal{R}_S^\sim, \Sigma_{1,S})$, where $\mathcal{R}_S^\sim \rightarrow S$ is a family of rubber expanded pairs, and $\Sigma_{1,S}$ is a section $S \rightarrow \mathcal{R}_S^\sim$ whose image in $\mathcal{R}_S^\sim$ is also denoted by $\Sigma_{1,S}$. Let $D_{0,S}, D_{\infty,S}$ be the relative divisors in $\mathcal{R}_S^\sim$. Note that $\Sigma_{1,S}$ is allowed to intersect the nodal and relative divisors.
In particular, on a geometric point $s \in S$, we have $\mathcal{R}_s^\sim \cong \Delta[k_s]_0$ for some $k_s \geq 0$, where
\[
\Delta[k_s]_0 = \Delta \cup_D \underbrace{\Delta \cup_D \cdots \cup_D \Delta}_{k_s\text{-times}}.
\]

The relative divisors in the rubber expanded pair are denoted by $D_{0,s}, D_{\infty,s}$, and $\Sigma_{1,s}$ is a point in $\Delta[k_s]_0$. We denote by $\Delta_{1,s}$ the first component $\Delta$ in $\Delta[k_s]_0$, which contains $D_{0,s}$, and denote by $D_{1,s}$ the nodal divisor where $\Delta_{1,s}$ meets other irreducible components.

\begin{defi}
Let $D(0|1,\infty)$ be the divisor in $\mathcal{R}^\sim$, whose set of $S$-points consists of pairs $(\mathcal{R}_S^\sim, \Sigma_{1,S})$ where $\Sigma_{1,s}$ does not lie in $\Delta_{1,s}\setminus D_{1,s}$ for any geometric point $s$ in $S$.

Let $R_-+R_+=R$ and we define $D(0|1,\infty)^{R_-,R_+}$ as the divisor in $\mcal R^{\sim,R}$ such that for a generic point in $D(0|1,\infty)$, the first $\Delta$ (that is, containing $D_0$) carries weight $R_-$ and the second $\Delta$ (that is, containing $\Sigma_{1,s}$ and $D_{\infty}$) carries weight $R_+$.
\end{defi}

\begin{nota}[{\cite[Rmk.~A.6]{CZZ24}}]
We let $\psi_0(\mcal R^{\sim})$ be the cotangent line bundle class that can be identified with the pullback of $\psi_0$
on $\mathfrak M_{0,3}^{\text{ss}}$ where we view $\mcal R^{\sim}$ as an open substack of $\mathfrak M^{\text{ss}}_{0,3}\times D$; see \cite[Rmk.~A.6]{CZZ24} for details.
\end{nota}

Then we have the following facts.
Let $(R_-,R_+)$ be a splitting datum of $R$, then we have $\psi_0(\mcal R^{\sim})=D(0|1,\infty)$
with decomposition of components
\[D(0|1,\infty)^R:=D(0|1,\infty)\times_{\mcal R^{\sim}}\mcal R^{\sim,R}
    =\bigcup_{R_-+R_+=R}D(0|1,\infty)^{R_-,R_+}\]
with $D(0|1,\infty)^{R_-,R_+}\cong\mcal A^{\sim,R_-}\times\mcal R^{\sim,\circ,R_+}$.
Which induce an isomorphism
\begin{equation}\label{eqn-CZZA10-4}
    \mcal Z^{\sim,R}\times_{\mcal R^{\sim,R}}D(0|1,\infty)^{R_-,R_+}\cong
    \Hilb^{\sim,R_-}(\Delta,D_0\sqcup D_{\infty})\times_{\spec\bb C}\mcal Z^{\sim,R_+}.
    \end{equation}
Now we have the following diagram
\begin{equation}
\begin{tiny}
\begin{tikzcd}
	{\mcal Z^{\sim,R}\times_{\mcal R^{\sim,R}}D(0|1,\infty)^{R_-,R_+}} & {\Hilb^{\sim,R}(\Delta,D_0\sqcup D_{\infty})\times_{\mcal A^{\sim,R}}D(0|1,\infty)^{R_-,R_+}} & {D(0|1,\infty)^{R_-,R_+}} \\
	{\mcal Z^{\sim,R}} & {\Hilb^{\sim,R}(\Delta,D_0\sqcup D_{\infty})\times_{\mcal A^{\sim,R}}\mcal R^{\sim,R}} & {\mcal R^{\sim,R}} \\
	& {\Hilb^{\sim,R}(\Delta,D_0\sqcup D_{\infty})} & {\mcal A^{\sim,R}} \\
	{\Hilb^{\sim,R_-}(\Delta,D_0\sqcup D_{\infty})\times\mcal Z^{\sim,R_+}} & {\Hilb^{\sim,R_-}(\Delta,D_0\sqcup D_{\infty})\times\Hilb^{\sim,R_+}(\Delta,D_0\sqcup D_{\infty})} & {\mcal A_0^{\sim,\dagger,R_-,R_+}.}
	\arrow[from=1-1, to=1-2]
	\arrow[from=1-1, to=2-1]
	\arrow["\Box"{description}, draw=none, from=1-1, to=2-2]
	\arrow["\cong"', bend right=30, from=1-1, to=4-1]
	\arrow[from=1-2, to=1-3]
	\arrow[hook, from=1-2, to=2-2]
	\arrow["\Box"{description}, draw=none, from=1-2, to=2-3]
	\arrow["{i_{D(0|1,\infty),{R_-,R_+}}}", hook, from=1-3, to=2-3]
	\arrow[from=2-1, to=2-2]
	\arrow["{\pi_{\mcal Z^{\sim}}}"', from=2-1, to=3-2]
	\arrow[from=2-2, to=2-3]
	\arrow[from=2-2, to=3-2]
	\arrow["\Box"{description}, draw=none, from=2-2, to=3-3]
	\arrow[from=2-3, to=3-3]
	\arrow[from=3-2, to=3-3]
	\arrow["\Box"{description}, draw=none, from=3-2, to=4-3]
	\arrow["{\text{id}\times\pi_{\mcal Z^{\sim}}}", from=4-1, to=4-2]
	\arrow[hook, from=4-2, to=3-2]
	\arrow[from=4-2, to=4-3]
	\arrow[hook, from=4-3, to=3-3]
\end{tikzcd}
\end{tiny}
\end{equation}
Then we have the following degeneration formula similar to those in the proof 
of \cite[Lem.~6.6]{CZZ24}:
\begin{lem}
We have:
\begin{enumerate}
    \item the identifications
    \begin{equation}\label{eqn-qdeuse-01}
        \begin{aligned}
            \psi_0(\mcal R^{\sim,R})\cap[\mcal Z^{\sim,R}]^{\vir}
            &=\sum_{R_-+R_+=R}i_{D(0|1,\infty),{R_-,R_+}}^![\mcal Z^{\sim,R}]^{\vir}\\
            &=\sum_{R_-+R_+=R}[\Hilb^{\sim,R_-}(\Delta,D_0\sqcup D_{\infty})]_{\bb T}^{\vir}\boxtimes[\mcal Z^{\sim,R_+}]^{\vir}.
        \end{aligned}
    \end{equation}
    \item If $\pi:\mcal R^{\sim}\to\mcal A^{\sim}$, then
    \begin{equation}\label{eqn-qdeuse-02}
        \psi_0(\mcal R^{\sim,R})\cap[\mcal Z^{\sim,R}]^{\vir}
        =\pi^*\psi_0(\mcal A^{\sim,R})\cap[\mcal Z^{\sim,R}]^{\vir}.
    \end{equation}
\end{enumerate}
\end{lem}
Using this degeneration formula, one can get the following result follows from a topological recursion relation
which is a quantum differential equation (\textsf{qde} for short).

\begin{prp}\label{prp-qde}
Define the rubber invariants
\[
W_{\infty}:=1+\sum_{n\geq1}\sum_{n_0+\cdots+n_{r-1}=n\leftrightarrow R}
q_0^{n_0}\cdots q_{r-1}^{n_{r-1}}\int_{[\Hilb^{\sim,R}(\Delta,D_0\sqcup D_{\infty})]_{\bb T}^{\vir}}
\frac{e_{\bb T\times\bb C_m^*}((\scr O^{[R]})^{\vee}\otimes e^m)}{s_3-\psi_0}
\]
and
\begin{align*}
F_{\infty,\ell}&:=\sum_{n\geq1}\sum_{n_0+\cdots+n_{r-1}=n\leftrightarrow R}
q_0^{n_0}\cdots q_{r-1}^{n_{r-1}}\int_{[\Hilb^{\sim,R}(\Delta,D_0\sqcup D_{\infty})]_{\bb T}^{\vir}}
e_{\bb T\times\bb C_m^*}((\scr O^{[R]})^{\vee}\otimes e^m)\psi_0^{\ell}\\
&\in\mathrm{\text{Frac}}\left(\frac{\bb C[s_1,s_2,s_3+s_4]}{s_1+s_2+s_3+s_4}\right)[m][[q_0,...,q_{r-1}]]
\cong\bb C(s_1,s_3+s_4)[m][[q_0,...,q_{r-1}]]
\end{align*}
for $\ell\geq0$, where for the rubber theories here, the bubble components $\Delta$ come from the chart of $X$ containing
$D_{\infty}$ and the torus on the fiber has equivariant weights $-s_1,-s_2,-(s_3+s_4)$.
Then we have
\begin{equation}
\log W_{\infty}=\frac{F_{\infty,0}}{s_3}.
\end{equation}
Here $\psi_0=\psi(\mcal A^{\sim})$ and we pull it back to $\Hilb^{\sim,R}(\Delta,D_0\sqcup D_{\infty})$
without changing the notation.
\end{prp}
\begin{proof}
We follow \cite[Lem.~4]{MNOP2}. In this proof, we use the notation
\[
Q\frac{\partial}{\partial Q}:=\sum_{i=0}^{r-1}q_i\frac{\partial}{\partial q_i}.
\]
Pulling $\mcal Z^{\sim,R}$ back along
$\bb C^2\to[\bb C^2/\bb Z_r]$ gives a finite flat family of degree $\dim R$.
Since this cover has degree $r$, the projection formula gives
\[
\pi_{\mcal Z^{\sim,R},*}[\mcal Z^{\sim,R}]^{\vir}
=\frac{\dim R}{r}[\Hilb^{\sim,R}(\Delta,D_0\sqcup D_{\infty})]_{\bb T}^{\vir}.
\]
For $\ell\geq1$ we have the following topological recursion \textsf{qde}:
\begin{smallalign}
\frac{1}{r}Q\frac{\partial}{\partial Q}F_{\infty,\ell}=
&\sum_{n\geq1}\sum_{n_0+\cdots+n_{r-1}=n\leftrightarrow R}
q_0^{n_0}\cdots q_{r-1}^{n_{r-1}}\int_{[\mcal Z^{\sim,R}]^{\vir}}
\pi_{\mcal Z^{\sim,R}}^*e_{\bb T\times\bb C_m^*}((\scr O^{[R]})^{\vee}\otimes e^m)\cdot\pi_{\mcal Z^{\sim,R}}^*\psi_0^{\ell}
\label{eqn-qde-01}\\
=&\sum_{n\geq1}\sum_{n_0+\cdots+n_{r-1}=n\leftrightarrow R}
q_0^{n_0}\cdots q_{r-1}^{n_{r-1}}\int_{[\mcal Z^{\sim,R}]^{\vir}}
\pi_{\mcal Z^{\sim,R}}^*e_{\bb T\times\bb C_m^*}((\scr O^{[R]}){^\vee}\otimes e^m)\cdot\pi_{\mcal Z^{\sim,R}}^*\psi_0^{\ell-1}\cdot\psi_0(\mcal R^{\sim})\label{eqn-qde-02}\\
=&\sum_{n\geq1}\sum_{n_0+\cdots+n_{r-1}=n\leftrightarrow R}
q_0^{n_0}\cdots q_{r-1}^{n_{r-1}}\sum_{R_-+R_+=R}\left(
\int_{[\Hilb^{\sim,R_-}(\Delta,D_0\sqcup D_{\infty})]_{\bb T}^{\vir}}
e_{\bb T\times\bb C_m^*}((\scr O^{[R_-]})^{\vee}\otimes e^m)\cdot\psi_0^{\ell-1}\right.\label{eqn-qde-03}\\
&\cdot\left.
\int_{[\mcal Z^{\sim,R_+}]^{\vir}}
\pi_{\mcal Z^{\sim,R_+}}^*e_{\bb T\times\bb C_m^*}((\scr O^{[R_+]})^{\vee}\otimes e^m)\right)
=\frac{F_{\infty,\ell-1}}{r}Q\frac{\partial}{\partial Q}F_{\infty,0},\nonumber
\end{smallalign}

\noindent where (\ref{eqn-qde-01}) follows from the projection formula above,
(\ref{eqn-qde-02}) from (\ref{eqn-qdeuse-02}), and (\ref{eqn-qde-03}) from (\ref{eqn-qdeuse-01}).

After cancelling $1/r$, the operator $Q\frac{\partial}{\partial Q}$ acts on a monomial
of total degree $n>0$ by multiplication by $n$. Since all $F_{\infty,\ell}$ have zero
constant term, induction on $\ell$ gives
\[
F_{\infty,\ell}=\frac{F_{\infty,0}^{\ell+1}}{(\ell+1)!},\quad\ell\geq0.
\]
Therefore,
\[
W_{\infty}=1+\sum_{\ell=0}^{\infty}\frac{F_{\infty,\ell}}{s_3^{\ell+1}}
=1+\sum_{\ell=0}^{\infty}\frac{F_{\infty,0}^{\ell+1}}{(\ell+1)!s_3^{\ell+1}}=\exp\left(\frac{F_{\infty,0}}{s_3}\right).
\]
Hence $\log W_{\infty}=\frac{F_{\infty,0}}{s_3}$.
\end{proof}

Next we need a localization formula.
Using the standard arguments as in \cite[Lem.~6.5]{CZZ24}, we have the following formula
through the torus localization by the sub-torus $\bb C^*=\{(1,1,t_3,t_3^{-1})\}\subset\bb T$.
\begin{prp}\label{prp-localization-formula}
Define the relative $\mathsf{DT}_4$-invariant of $(X,D_{\infty})$ as
\begin{align*}
    &\DTC(X,D_{\infty}):=1+\sum_{n\geq1}\sum_{n_0+\cdots+n_{r-1}=n\leftrightarrow R}
q_0^{n_0}\cdots q_{r-1}^{n_{r-1}}\int_{[\Hilb^{R}(X,D_{\infty})]_{\bb T}^{\vir}}
e_{\bb T\times\bb C_m^*}((\scr O^{[R]})^{\vee}\otimes e^m),
\end{align*}
then we have 
\begin{equation}
\DTC(X,D_{\infty})=W_{\infty}\cdot\DTC_{[\bb C^4/\bb Z_r],\scr O,o}.
\end{equation}
\end{prp}

Finally, we can use the results of this section to prove the following result.

\begin{prp}\label{prp-pole-s3-logDT4}
$\log\DTC_{[\bb C^4/\bb Z_r],\scr O,o}$ has a pole of order $1$ along $s_3=0$.
\end{prp}
\begin{proof}
By Proposition \ref{prp-qde} and Proposition \ref{prp-localization-formula}, we have
\[\log\DTC_{[\bb C^4/\bb Z_r],\scr O,o}
=\log\DTC(X,D_{\infty})-\frac{F_{\infty,0}}{s_3}.\]
Pulling a relative subscheme back along
$\bb C^2\to[\bb C^2/\bb Z_r]$ gives a $\bb Z_r$-invariant relative subscheme
of length $\dim R$ on $\bb C^2\times\Tot(\scr O_{\bb P^1}(-1))$.
Collapsing bubbles, applying Hilbert--Chow and blowing down give proper
$\bb T$-equivariant maps
\[
\Hilb^R(X,D_{\infty})
\to\Sym^{\dim R}\big(\bb C^2\times\Tot(\scr O_{\bb P^1}(-1))\big)
\to\Sym^{\dim R}(\bb C^4).
\]
Here the first map is the relative Hilbert--Chow morphism restricted to the closed
locus of $\bb Z_r$-invariant subschemes of representation type $R$.
The tangent weights on $\bb C^4$ are $-s_1,-s_2,-s_4,-(s_3+s_4)$.
As in \cite[Lem.~3]{MNOP2}, compose with the finite map
\[
\begin{aligned}
\Sym^{\dim R}(\bb C^4)&\longrightarrow\bigoplus_{\ell=1}^{\dim R}\bb C^4,\\
\{(x_i,y_i,z_i,w_i)\}&\longmapsto
\bigoplus_{\ell=1}^{\dim R}\left(\sum_ix_i^{\ell},\sum_iy_i^{\ell},\sum_iz_i^{\ell},\sum_iw_i^{\ell}\right).
\end{aligned}
\]
Finiteness follows from Newton identities, applied to each coordinate.
Giving the $\ell$-th summand weights
$-\ell s_1,-\ell s_2,-\ell s_4,-\ell(s_3+s_4)$ makes this map $\bb T$-equivariant.

Pushing the virtual class with its tautological insertion to this representation shows
that the coefficients of $\DTC(X,D_{\infty})$ and its formal logarithm 
only have poles at $s_1,s_2,s_4,s_3+s_4$.
Since $s_3+s_4=-(s_1+s_2)$ and
$F_{\infty,0}\in\bb C(s_1,s_2)[m][[q_0,...,q_{r-1}]]$,
the polar part along $s_3=0$ is $-F_{\infty,0}/s_3$.
For the single-box partition $\pi=\{(0,0,0,0)\}$, Corollary \ref{coro-sign-rule-oB} gives
that the coefficient of $q_0$ in 
$\log\mathsf{DT}_{[\mathbb C^4/\mathbb Z_r],\mathscr O,o}$ is
$-\frac{m}{s_3}-\frac{m}{s_4}$.
Then the coefficient of $q_0$ in 
$F_{\infty,0}$ is $m\neq0$. Hence $F_{\infty,0}\neq0$ and the result follows.
\end{proof}

\subsection{Two compactifications}\label{subsection-toric-cpts}
Motivated by \cite{Zhou18-2} and the Donaldson-Thomas crepant resolution conjecture stated in \cite[Conj.~5.16,~6.7]{CKM23}, we need to compare the invariants of the following two crepant resolutions
\[
[\mathbb{C}^4/\mathbb{Z}_r] \longrightarrow \mathbb{C}^4/_{\text{aff}}\mathbb{Z}_r \longleftarrow \widetilde{\mathbb{C}^2/\mathbb{Z}_r}\times\mathbb{C}^2
\]
of the affine GIT quotient $\mathbb{C}^4/_{\text{aff}}\mathbb{Z}_r$ (which is an 
$A_{r-1}$-type singularity). The two resolutions are isomorphic outside the preimages of 
the singular locus $\{0\}\times\bb C^2$; the exceptional chain and the stacky point
are both taken in product with $\bb C^2$.

We first construct proper toric compactifications of $\widetilde{\bb C^2/\bb Z_r}$
and $[\bb C^2/\bb Z_r]$, and then take their canonical-bundle total spaces in product
with $\bb C$ to obtain the required partial compactifications.
Consider $X_0:=\bb P^1\times\bb P^1$ with divisors $V_0=\{0\}\times\bb P^1$,
$H_0=\bb P^1\times\{0\}$, and $H_{\infty}=\bb P^1\times\{\infty\}$.

Let $q_1=(\infty,0)\in X_0$ and consider $X_1:=\mathrm{Bl}_{q_1}(X_0)$ with exceptional divisor 
$E_1$. Then we let $q_2:=E_1\cap H_0$ and consider $X_2:=\mathrm{Bl}_{q_2}(X_1)$.
We can repeat this process to obtain $X_r$ with exceptional divisors $E_1,\dots,E_r$,
where for $i\geq2$ we blow up the intersection point $E_{i-1}\cap H_0$.
Here we still use $H_0$, etc. to denote their strict transforms.
Below is a schematic diagram of $X_r$ and its corresponding fan:
\begin{figure}[!ht]
    \centering
    \tikzset{every picture/.style={line width=0.6pt}} 
\begin{tikzpicture}[x=0.6pt,y=0.6pt,yscale=-1,xscale=1]

\draw    (461.5,202) -- (538.5,202) ;
\draw [shift={(540.5,202)}, rotate = 180] [color={rgb, 255:red, 0; green, 0; blue, 0 }  ][line width=0.75]    (10.93,-3.29) .. controls (6.95,-1.4) and (3.31,-0.3) .. (0,0) .. controls (3.31,0.3) and (6.95,1.4) .. (10.93,3.29)   ;
\draw    (461.5,202) -- (461.5,279) ;
\draw [shift={(461.5,281)}, rotate = 270] [color={rgb, 255:red, 0; green, 0; blue, 0 }  ][line width=0.75]    (10.93,-3.29) .. controls (6.95,-1.4) and (3.31,-0.3) .. (0,0) .. controls (3.31,0.3) and (6.95,1.4) .. (10.93,3.29)   ;
\draw    (461.5,202) -- (383.5,202) ;
\draw [shift={(381.5,202)}, rotate = 360] [color={rgb, 255:red, 0; green, 0; blue, 0 }  ][line width=0.75]    (10.93,-3.29) .. controls (6.95,-1.4) and (3.31,-0.3) .. (0,0) .. controls (3.31,0.3) and (6.95,1.4) .. (10.93,3.29)   ;
\draw    (461.5,202) -- (461.5,124) ;
\draw [shift={(461.5,122)}, rotate = 90] [color={rgb, 255:red, 0; green, 0; blue, 0 }  ][line width=0.75]    (10.93,-3.29) .. controls (6.95,-1.4) and (3.31,-0.3) .. (0,0) .. controls (3.31,0.3) and (6.95,1.4) .. (10.93,3.29)   ;
\draw    (461.5,202) -- (538.09,125.41) ;
\draw [shift={(539.5,124)}, rotate = 135] [color={rgb, 255:red, 0; green, 0; blue, 0 }  ][line width=0.75]    (10.93,-3.29) .. controls (6.95,-1.4) and (3.31,-0.3) .. (0,0) .. controls (3.31,0.3) and (6.95,1.4) .. (10.93,3.29)   ;
\draw    (461.5,202) -- (541.72,13.84) ;
\draw [shift={(542.5,12)}, rotate = 113.09] [color={rgb, 255:red, 0; green, 0; blue, 0 }  ][line width=0.75]    (10.93,-3.29) .. controls (6.95,-1.4) and (3.31,-0.3) .. (0,0) .. controls (3.31,0.3) and (6.95,1.4) .. (10.93,3.29)   ;
\draw   (87.43,56) -- (303.58,56) -- (303.58,264) -- (87.43,264) -- cycle ;
\draw    (125.2,56) -- (125.2,264) ;
\draw    (88.05,231.81) -- (303.48,231.81) ;
\draw    (86.81,83.24) -- (303.48,83.24) ;
\draw    (169.77,65.9) -- (216.81,119.14) ;
\draw    (209.39,96.86) -- (179.67,129.05) ;
\draw    (175.96,111.71) -- (206.91,147.62) ;
\draw    (175.96,156.29) -- (201.96,183.52) ;
\draw    (199.48,168.67) -- (169.77,204.57) ;
\draw  [dash pattern={on 4.5pt off 4.5pt}]  (167.29,187.24) -- (204.43,262.76) ;

\draw (520,80.4) node [anchor=north west][inner sep=0.75pt]    {$\vdots $};
\draw (437,130.4) node [anchor=north west][inner sep=0.75pt]    {$H_{0}$};
\draw (431,247.4) node [anchor=north west][inner sep=0.75pt]    {$H_{\infty }$};
\draw (391,176.4) node [anchor=north west][inner sep=0.75pt]    {$V_{0}$};
\draw (513,207.4) node [anchor=north west][inner sep=0.75pt]    {$V_{\infty }$};
\draw (530,141.4) node [anchor=north west][inner sep=0.75pt]    {$E_{1}$};
\draw (539,24.4) node [anchor=north west][inner sep=0.75pt]    {$E_{r}$};
\draw (151,59.4) node [anchor=north west][inner sep=0.75pt]    {$E_{r}$};
\draw (227,62.4) node [anchor=north west][inner sep=0.75pt]    {$H_{0}$};
\draw (222,212.4) node [anchor=north west][inner sep=0.75pt]    {$H_{\infty }$};
\draw (101,148.4) node [anchor=north west][inner sep=0.75pt]    {$V_{0}$};
\draw (167,237.4) node [anchor=north west][inner sep=0.75pt]    {$V_{\infty }$};
\draw (38,55.4) node [anchor=north west][inner sep=0.75pt]    {$X_{r} :$};
\draw (167,92.4) node [anchor=north west][inner sep=0.75pt]    {$E_{r-1}$};
\draw (210,135.4) node [anchor=north west][inner sep=0.75pt]    {$E_{r-2}$};
\draw (156,143.4) node [anchor=north west][inner sep=0.75pt]    {$E_{2}$};
\draw (181,185.4) node [anchor=north west][inner sep=0.75pt]    {$E_{1}$};
\end{tikzpicture}
\caption{The schematic diagram and fan of $X_r$.}
\label{Fig-Xr-fan}
\end{figure}
Here $V_{\infty}$ is the strict transform of $\{\infty\}\times\bb P^1$.

Let $U_r:=X_r\backslash(V_0\cup H_0\cup H_{\infty})$. By the fan of $X_r$
in Figure \ref{Fig-Xr-fan}, one has $U_r\cong\widetilde{\bb C^2/\bb Z_r}$.

To construct $Y_r$, contract the chain of $(-2)$-curves $E_1\cup\cdots\cup E_{r-1}$
and take the canonical smooth Deligne-Mumford stack of the resulting projective toric
surface; see \cite{FMM10}. In terms of Figure \ref{Fig-Xr-fan}, this amounts to removing
the rays corresponding to $E_1,\dots,E_{r-1}$. The cone between the rays corresponding
to $V_{\infty}$ and $E_r$ gives the open substack $U_r'\cong[\bb C^2/\bb Z_r]$.
The contraction is crepant and is an isomorphism away from the exceptional chain.
The four maximal toric charts outside $U_r$ and $U_r'$ are identical.

We will use the absolute $\mathsf{DT}_4$-invariants of
\[
\Tot(\omega_{X_r})\times\bb C\quad\text{and}\quad\Tot(\omega_{Y_r})\times\bb C.
\]
Since $\omega_{U_r}$ and $\omega_{U_r'}$ are trivial, these spaces contain
$\widetilde{\bb C^2/\bb Z_r}\times\bb C^2$ and $[\bb C^4/\bb Z_r]$ as open substacks.
Lift the surface torus to the canonical bundle, scale its fibre by $\lambda$,
and scale the last coordinate by $\lambda^{-1}$. Using the covering torus on
$[\bb C^2/\bb Z_r]$, the characters on the open chart are
\[
(t_1,t_2,t_3,t_4)=(t_1,t_2,(t_1t_2)^{-1}\lambda,\lambda^{-1}).
\]
This gives the original Calabi-Yau torus $\bb T$ with $s_1+s_2+s_3+s_4=0$.
We choose the orientations from Lemma \ref{lem-ori-identification}, viewing these
fourfolds as canonical-bundle total spaces over $\Tot(\omega_{X_r})$ and
$\Tot(\omega_{Y_r})$, respectively, with the last coordinate as the canonical fibre.
These orientations restrict to the previously chosen orientations on the original
open charts and agree on the common charts.

\begin{lem}\label{lem-comparison-poles}
$\log\DTC_{[\bb C^4/\bb Z_r],\scr O,o}-\log\DTC_{\widetilde{\bb C^2/\bb Z_r}\times\bb C^2,\scr O}$
only has poles along $s_3=0$ and $s_4=0$.
\end{lem}
\begin{proof}
We identify $q=Q=q_0\cdots q_{r-1}$ on the smooth side, as in Remark
\ref{Rmk-Change-variables}. Torus localization expresses the absolute invariants
of the two canonical-bundle total spaces as products of vertex contributions.
The obstruction complexes and tautological bundles split over the disjoint supports,
and the orientations are compatible with this splitting by the canonical-bundle
presentations above, as in \S\ref{sect on DT toric 4folds}.
The four common smooth charts have identical equivariant geometry, orientations,
and insertions, and their counting variable is $Q$ on both sides. Cancelling their
vertex contributions and taking logarithms gives
\begin{align*}
&\log\DTC_{[\bb C^4/\bb Z_r],\scr O,o}
-\log\DTC_{\widetilde{\bb C^2/\bb Z_r}\times\bb C^2,\scr O}\\
&=\log\DTC_{\Tot(\omega_{Y_r})\times\bb C,\scr O}
-\log\DTC_{\Tot(\omega_{X_r})\times\bb C,\scr O}.
\end{align*}

Consider the sub-torus $\bb C^*=\{(1,1,t_3,t_3^{-1})\}\subset\bb T$.
Its fixed loci in the two total spaces are the proper zero sections $X_r$ and $Y_r$.
The fixed loci in the Hilbert spaces of any fixed numerical class are also proper.
Indeed, a $\bb C^*$-invariant zero-dimensional subscheme has finite invariant support,
so its support lies in the zero section. For a fixed class, its length is bounded,
and it is contained in a bounded infinitesimal neighbourhood of that section,
which is proper. On the stack side, the same argument uses a generating bundle
and the properness of the corresponding Quot space; see \cite{OS03}.

Localization by this sub-torus, retaining $\bb T\times\bb C_m^*$-equivariance,
only requires inverting characters whose restrictions to $\bb C^*$ are nonzero.
Since the fixed loci are proper, the coefficients of both absolute invariants
and their formal logarithms belong to
\[
\bb C[s_1,s_2,s_3,m]
\left[\frac{1}{a s_1+b s_2+c s_3}:a,b,c\in\bb Q,\ c\neq0\right],
\]
where $s_4=-(s_1+s_2+s_3)$. Thus no nonzero linear form in $s_1,s_2$ occurs
in the reduced denominator of any coefficient of the difference.

On the other hand, the proof of Proposition \ref{prp-pole-s3-logDT4} gives
\[
\log\DTC_{[\bb C^4/\bb Z_r],\scr O,o}
=\log\DTC(X,D_{\infty})-\frac{F_{\infty,0}}{s_3},
\]
where the coefficients of $\DTC(X,D_{\infty})$ only have poles in monomials of
$s_1,s_2,s_4,s_3+s_4$. Since $s_3+s_4=-(s_1+s_2)$ and
$F_{\infty,0}\in\bb C(s_1,s_2)[m][[q_0,...,q_{r-1}]]$, and localization in the
rubber theory only inverts nonzero characters of its two-dimensional torus,
the denominator of each coefficient of $\log\DTC_{[\bb C^4/\bb Z_r],\scr O,o}$
is a product of powers of $s_3,s_4$ and nonzero linear forms in $s_1,s_2$.
By Corollary \ref{coro-DT-ArC2}, the same is true of
$\log\DTC_{\widetilde{\bb C^2/\bb Z_r}\times\bb C^2,\scr O}$.

Then the result follows from combining these two observations.
\end{proof}

\begin{nota}\label{Rmk-Change-variables}
When comparing the invariants $\DTC_{[\bb C^4/\bb Z_r],\scr O,o}$ and 
$\DTC_{\widetilde{\bb C^2/\bb Z_r}\times\bb C^2,\scr O}$, a compatible change of 
variables is given by $q\mapsto Q=q_0\cdots q_{r-1}$. This relation arises because 
the non-commutative models for $\Hilb^n\left(\widetilde{\bb C^2/\bb Z_r}\times\bb C^2\right)$ and 
$\Hilb^R([\bb C^4/\bb Z_r])$ are described by the same framed quiver in Figure \ref{fig-framed-quiver}, 
but with different dimension vectors: $(n, n, \dots, n, 1)$ for the former, and 
$(n_0, n_1, \dots, n_{r-1}, 1)$ for the latter.
Geometrically, a point in the non-stacky locus lifts to a free $\bb Z_r$-orbit,
whose structure sheaf carries the regular representation
$\rho_0\oplus\cdots\oplus\rho_{r-1}$. Hence contributions from the common smooth
charts are series in $Q$ on both sides.
\end{nota}

\subsection{\texorpdfstring{$\mathsf{DT}_4$}{DT4}-invariants of \texorpdfstring{$[\mathbb C^4/\bb Z_r]$}{[C4/Zr]}}\label{sec-DT4-C4Zr-final}
In this section, we will prove our main theorem using the results in the previous sections. Firstly, we need the following lemma.

\begin{lem}\label{lem-divide-property}
The coefficients of $q_0^{n_0}\cdots q_{r-1}^{n_{r-1}}$ in $\DTC_{[\bb C^4/\bb Z_r],\scr O,o}$ are divided by $m$ and $s_1+s_2$
for $n_0+\cdots +n_{r-1}>0$.
\end{lem}
\begin{proof}
We will use a similar idea to \cite[Lem.~5]{MNOP2}
via cohomological limits.

For $m$, note that the only term that contains $y$ is 
\[
-y\overline{Z_{\pi}^{(0)}}
=-y\sum_{\substack{(i,j,k,l)\in\pi\\i-j\equiv0\mmod r}}
t_1^{-i}t_2^{-j}t_3^{-k}t_4^{-l}.
\]
Since $(0,0,0,0)\in\pi$ for $\pi\neq\emptyset$, the tautological
insertion contains a factor $m$. The denominator is independent of $m$,
so each vertex contribution is divisible by $m$.

For $s_1+s_2$, we only need to prove the negativity of the sum of the coefficients of $t_1^it_2^it_3^0 t_4^0y^0$
in the K-theoretic invariants.
This is the constant term of $\widetilde{v}^{\DT,\bb Z_r}_{\pi}(x,x^{-1},t_3,t_3^{-1},y)$ for any solid partition $\pi\neq\emptyset$. By direct calculation, the constant term of $\widetilde{v}^{\DT,\bb Z_r}_{\pi}(x,x^{-1},t_3,t_3^{-1},y)$ is 
\begin{align*}
    C^{\pi}=&a_{00}^{\pi}-\sum_{m,k\in\bb Z}\left(2(a_{km}^{\pi})^2-a_{k+1,m}^{\pi}a_{km}^{\pi}-a_{k-1,m}^{\pi}a_{km}^{\pi}\right.\\
    &\left.-2
    a_{km}^{\pi}a_{k,m-1}^{\pi}+a_{k+1,m}^{\pi}a_{k,m-1}^{\pi}+
   a_{k-1,m}^{\pi}a_{k,m-1}^{\pi} \right)\\
   =&a_{00}^{\pi}-\frac{1}{2}\sum_{m,k\in\bb Z}\left(a_{km}^{\pi}-a_{k+1,m}^{\pi}-a_{k,m+1}^{\pi}+a_{k+1,m+1}^{\pi}\right)^2,
\end{align*}
where $a_{ij}^{\pi}:=|\{(a+i,a,b+j,b)\in\pi\}|$.
Since all these are integers, we have the inequality
\begin{align*}
   C^{\pi}\leq a_{00}^{\pi}-\frac{1}{2}\sum_{m,k\in\bb Z}\left|a_{km}^{\pi}-a_{k+1,m}^{\pi}-a_{k,m+1}^{\pi}+a_{k+1,m+1}^{\pi}\right|.
\end{align*}
By splitting them into four quadrants with different signs, we have
\begin{align*}
   C^{\pi}\leq& a_{00}^{\pi}-\frac{1}{2}\sum_{m\geq0,k\geq0}(a_{km}^{\pi}-a_{k+1,m}^{\pi}-a_{k,m+1}^{\pi}+a_{k+1,m+1}^{\pi})\\
   &-
   \frac{1}{2}\sum_{m<0,k<0}(a_{km}^{\pi}-a_{k+1,m}^{\pi}-a_{k,m+1}^{\pi}+a_{k+1,m+1}^{\pi})\\
   &+\frac{1}{2}\sum_{m\geq0,k<0}(a_{km}^{\pi}-a_{k+1,m}^{\pi}-a_{k,m+1}^{\pi}+a_{k+1,m+1}^{\pi})\\
   &+
   \frac{1}{2}\sum_{m<0,k\geq0}(a_{km}^{\pi}-a_{k+1,m}^{\pi}-a_{k,m+1}^{\pi}+a_{k+1,m+1}^{\pi})\\
   =&a_{00}^{\pi}-\frac{4a_{00}^{\pi}}{2}=-a_{00}^{\pi}<0,
\end{align*}
since $|\pi|>0$. This proves the claim.
\end{proof}

Finally, combining all results established in previous sections yields the following main theorem of this paper, which confirms Conjecture \ref{conj-DT-C4Zr}.

\begin{thm}\label{thm-main-DT4-C4Zr}
We have
\begin{equation}
\begin{aligned}
    \DTC_{[\bb C^4/\bb Z_r],\scr O,o}(m,q_0,...,q_{r-1})
    =&M(1,-Q)^{-\frac{m}{s_4}\left(\frac{r(s_1+s_2)}{s_3}+\frac{(s_1+s_2)(s_1+s_2+s_3)}{rs_1s_2}\right)}\\
    &\cdot\prod_{0<i\leq j<r}\widetilde{M}(q_{[i,j]},- Q)^{-\frac{m(s_1+s_2)}{s_3s_4}}.
\end{aligned}
\end{equation}
\end{thm}
\begin{proof}
If we let $\deg m=\deg s_i=1$ and $\deg q=0$, then the invariants
$\DTC_{[\bb C^4/\bb Z_r],\scr O,o}$ and $\DTC_{\widetilde{\bb C^2/\bb Z_r}\times\bb C^2,\scr O}$ are homogeneous of degree $0$
due to vertex formalisms.
Therefore, combining Proposition \ref{prp-pole-s3-logDT4}, Lemma \ref{lem-comparison-poles} and Lemma \ref{lem-divide-property}, we have
\begin{equation}\label{eqn-frac-DTDT}
\frac{\DTC_{[\bb C^4/\bb Z_r],\scr O,o}}{\DTC_{\widetilde{\bb C^2/\bb Z_r}\times\bb C^2,\scr O}}=F(q_0,...,q_{r-1})^{\frac{m(s_1+s_2)}{s_3s_4}}
\end{equation}
for some formal power series
$F\in1+(q_0,...,q_{r-1})\bb C[[q_0,...,q_{r-1}]]$
independent of $m$ and $s_i$. By Corollary \ref{coro-DT-ArC2} and Eqn. (\ref{eqn-frac-DTDT}), we have
\begin{equation}\label{eqn-DT-C4Zr-AC2-comp}
\begin{aligned}
\DTC_{[\bb C^4/\bb Z_r],\scr O,o}
&=F^{\frac{m(s_1+s_2)}{s_3s_4}}\cdot\DTC_{\widetilde{\bb C^2/\bb Z_r}\times\bb C^2,\scr O}\\
&=F^{\frac{m(s_1+s_2)}{s_3s_4}}\cdot M(-Q)^{-\frac{m}{s_4}
    \left(\frac{r(s_1+s_2)}{s_3}+\frac{(s_1+s_2)(s_1+s_2+s_3)}{rs_1s_2}\right)}.
\end{aligned}
\end{equation}
Taking $m=s_4$, we have
\begin{align*}
    \DTC_{[\bb C^4/\bb Z_r],\scr O,o}|_{m=s_4}=&\sum_{\pi}\left.\left((-1)^{\mu_{\pi}}
    e_{\bb T\times\bb C_m^*}\left(-\widetilde{v}_{\pi}^{\DT,\bb Z_r}\right)\right)\right|_{m=s_4}q_0^{|\pi|_{\rho_0}}\cdots q_{r-1}^{|\pi|_{\rho_{r-1}}}\\
   =&\sum_{\pi^{\text{3D}}}
    e_{\bb T}\left(-V_{\pi^{\text{3D}}}^{\DT,\bb Z_r}\right)q_0^{|\pi^{\text{3D}}|_{\rho_0}}\cdots q_{r-1}^{|\pi^{\text{3D}}|_{\rho_{r-1}}}\\
    =&\DTC_{[\bb C^3/\bb Z_r]}(q_0,...,q_{r-1}),
\end{align*}
where the second equality follows from the dimension reduction formula in \cite[Prop.~2.1]{CKM22} and the sign rules in Corollary \ref{coro-sign-rule-oB}, since $\mu_{\pi^{\text{3D}}}=0$ for any plane partition $\pi^{\text{3D}}$.
Hence, by \cite[Thm.~5.3]{Zhou18-2}, we have 
\begin{equation}\label{eqn-DT-C4Zr-m=s4}
\begin{aligned}
    &\DTC_{[\bb C^4/\bb Z_r],\scr O,o}|_{m=s_4}=\DTC_{[\bb C^3/\bb Z_r]}(q_0,...,q_{r-1})\\
    =&M(1,-Q)^{-\frac{r(s_1+s_2)}{s_3}-\frac{(s_1+s_2)(s_1+s_2+s_3)}{rs_1s_2}}
    \cdot\prod_{0<i\leq j<r}\widetilde{M}(q_{[i,j]},-Q)^{-\frac{s_1+s_2}{s_3}}.
\end{aligned}
\end{equation}
Therefore, we can take $m=s_4$ in Eqn. (\ref{eqn-DT-C4Zr-AC2-comp}) and use Eqn. (\ref{eqn-DT-C4Zr-m=s4}) to deduce that
\[F=\prod_{0<i\leq j<r}\widetilde{M}(q_{[i,j]},-Q)^{-1}.\] 
This gives the result.
\end{proof}

\begin{nota}\label{rmk-DT-C4Zr-different-sign}
Note that in the original conjecture in \cite[Cor.~6.6]{CKM23}, they use the sign rule $|\pi|_{R_0}+\mu_{\pi}$. This can be done if we consider $\scr O^{[R]}\otimes y^{-1}$ instead of $(\scr O^{[R]})^{\vee}\otimes y$ and consider the same vertices, since $\rank\scr O^{[R]}|_{p_{\pi}}=|\pi|_{R_0}$.
\end{nota}

As corollaries of our main theorem, we obtain the following two results.
\begin{coro}
Define the following insertion-free invariant:
\[
\DT_{[\bb C^4/\bb Z_r]}^{\mathsf{free}}:=1+\sum_{n_0+\cdots+n_{r-1}=n\leftrightarrow R,n>0}
q_0^{n_0}\cdots q_{r-1}^{n_{r-1}}\int_{[\Hilb^{R}([\bb C^4/\bb Z_r])]_{\bb T}^{\vir}}
1.
\]
Then we have
\begin{align*}
    \DT_{[\bb C^4/\bb Z_r]}^{\mathsf{free}}=&\exp\left(\frac{Q}{s_4}\left(\frac{r(s_1+s_2)}{s_3}+\frac{(s_1+s_2)(s_1+s_2+s_3)}{rs_1s_2}\right)\right)\\
    &\cdot\prod_{0<i\leq j<r}\exp\left((q_{[i,j]}+q_{[i,j]}^{-1})\frac{Q(s_1+s_2)}{s_3s_4}\right).
\end{align*}
\end{coro}
\begin{proof}
This proves \cite[Cor.~6.10]{CKM23}, see also \cite[Prop.~4.27]{ST23} as the partition functions of the pure gauge theories on the quotient stacks.
Using the insertion-free limit as in \cite[\S~6.3]{CKM23}, we have 
\[
\DT_{[\bb C^4/\bb Z_r]}^{\mathsf{free}}=\lim_{m\to\infty}\DT_{[\bb C^4/\bb Z_r]}(m,q_0',q_1,...,q_{r-1})|_{q_0=q_0'm}.
\]
Then the result follows from Theorem \ref{thm-main-DT4-C4Zr} and a direct computation.
\end{proof}

\begin{coro}\label{coro-Winfty-relative-1}
With the notations introduced in \S~\ref{sec_degeneration_formula}, we have 
\[
W_{\infty}=\left(M(1,-Q)^r
    \cdot\prod_{0<i\leq j<r}\widetilde{M}(q_{[i,j]},- Q)\right)^{-\frac{m}{s_3}}
\]
and
\[
\DT(X,D_{\infty})=M(1,-Q)^{\frac{m}{s_4}\left(r-\frac{(s_1+s_2)(s_1+s_2+s_3)}{rs_1s_2}\right)}
    \cdot\prod_{0<i\leq j<r}\widetilde{M}(q_{[i,j]},- Q)^{\frac{m}{s_4}}. 
\]
\end{coro}
\begin{proof}
By the arguments in \S~\ref{sec_degeneration_formula}, we have
\[
F_{\infty,0}=-\mathrm{Res}_{s_3=0}\left(\log\DT_{[\bb C^4/\bb Z_r]}\,ds_3\right),
\]
where the residue is taken at $s_3=0$, with $s_1,s_2$ fixed and $s_4=-(s_1+s_2+s_3)$.
On the other hand, by Theorem \ref{thm-main-DT4-C4Zr}, we have 
\[
\mathrm{Res}_{s_3=0}\left(\log\DT_{[\bb C^4/\bb Z_r]}\,ds_3\right)
=m\left(r\log M(1,-Q)+\sum_{0<i\leq j<r}\log\widetilde{M}(q_{[i,j]},- Q)\right).
\]
Then the expression for $W_{\infty}$ follows from these equations together with the relation
$W_{\infty}=\exp\left(\frac{F_{\infty,0}}{s_3}\right)$. Finally the second identity follows from the equality
$\DT(X,D_{\infty})=W_{\infty}\cdot \DT_{[\bb C^4/\bb Z_r]}$.
\end{proof}

Finally, as in \cite[Cor.~6.11]{CZZ24}, we can compute further invariants 
of the pair $(X,D_{\infty})$ with other insertions.
Indeed, let $D_{0}=[\bb C^2/\bb Z_r]\times\{0\}\times\bb C\subset X$. Consider a 
$(\bb C^*)^4\times\bb C^*_m$-equivariant line bundle
\[
\scr O[\ell]\in\mathrm{Pic}^{(\bb C^*)^4\times\bb C^*_m}(X)
\]
such that $\scr O[\ell]|_{X\backslash D_{\infty}}=\scr O$
and $\scr O[\ell]|_{X\backslash D_{0}}=\scr O\cdot t_3^{\ell}$.

\begin{coro}
Define the following invariant with a non-trivial tautological insertion:
\[
\DT(X,D_{\infty})[\ell]:=1+\sum_{n_0+\cdots+n_{r-1}=n\leftrightarrow R}
q_0^{n_0}\cdots q_{r-1}^{n_{r-1}}\int_{[\Hilb^{R}(X,D_{\infty})]_{\bb T}^{\vir}}
e_{\bb T\times\bb C_m^*}((\scr O[\ell]^{[R]})^{\vee}\otimes e^m).
\]
Then we have 
\[
\DT(X,D_{\infty})[\ell]=\DT(X,D_{\infty})\cdot\left(M(1,-Q)^r
    \cdot\prod_{0<i\leq j<r}\widetilde{M}(q_{[i,j]},- Q)\right)^{\ell}.
\]
\end{coro}
\begin{proof}
Indeed, by arguments similar to those in Proposition \ref{prp-localization-formula}, we obtain 
\[
\DT(X,D_{\infty})[\ell]=\DT_{[\bb C^4/\bb Z_r]}\cdot W_{\infty}|_{m\mapsto m-\ell s_3}.
\]
Then the result follows from this and Corollary \ref{coro-Winfty-relative-1}.
\end{proof}

\end{document}